\documentclass[12pt]{amsart}
\usepackage[frenchb]{babel}
\usepackage{a4wide}

\newcommand{\al}{\alpha}
\newcommand{\io}{\infty}
\newcommand{\dd}{\textup{d}}

\newcommand{\be}{\beta}
\newcommand{\ep}{\epsilon}

\newcommand{\de}{\delta}

\newcommand{\dis}{\displaystyle}
\def\poq#1#2{(#1;q)_#2}

\newtheorem{theo}{Th{\'e}or{\`e}me}
\newtheorem{lem}{Lemme}
\newtheorem{prop}{Proposition}

\theoremstyle{remark}

\newtheorem*{remark}{Remarque}
\newtheorem*{remarks}{Remarques}

\numberwithin{equation}{section}

\begin{document}

\title[]{
S{\'e}ries hyperg{\'e}om{\'e}triques basiques,
$q$-analogues des valeurs de la fonction z{\^e}ta et s{\'e}ries d'Eisenstein}
\author[]{C. Krattenthaler$^\dagger$, T. Rivoal et W. Zudilin$^\ddagger$}

\address{
Institut Girard Desargues,
Universit{\'e} Claude Bernard Lyon-I,
21, avenue Claude Bernard,
F-69622 Villeurbanne Cedex, France}
\email{kratt@igd.univ-lyon1.fr}
\address{
Laboratoire de Math{\'e}matiques Nicolas Oresme, CNRS UMR 6139,
Universit{\'e} de Caen,  BP 5186,
14032 Caen cedex,
France}
\email{rivoal@math.unicaen.fr}
\address{
Universit{\'e} d'{\'E}tat Lomonosov, Moscou,
Facult{\'e} de M{\'e}canique et Math{\'e}matiques,
Vorobiovy Gory, GSP-2, Moscou 119992 Russie}
\email{wadim@ips.ras.ru}
\thanks{$^\dagger$ Recherche partiellement support{\'e}e par la Fondation
Autrichienne de la Recherche Scientifique FWF, contrat P12094-MAT,
et par le Programme {\og Accro{\^\i}tre le potentiel humain de recherche \fg}
de la Commission Europ{\'e}enne, contrat HPRN-CT-2001-00272,
``Algebraic Combinatorics in Europe"}
\thanks{$^\ddagger$ Recherche support{\'e}e par une bourse {\og Alexander
von Humboldt \fg} et partiellement support{\'e}e par la Fondation Russe
pour la Recherche Fondamentale, contrat 03-01-00359}

\date{}



\begin{abstract}
Nous {\'e}tudions la nature arithm{\'e}tique    
de $q$-analogues des valeurs $\zeta(s)$ de la fonction
z{\^e}ta de Riemann, 
notamment des valeurs des fonctions $\zeta_q(s)=
\sum_{k=1} ^{\infty}q^k
\sum_{d\mid k} ^{}d^{s-1}$, $s=1,2,\dots$, o{\`u} $q$ est un nombre
complexe, $\vert q\vert<1$ (ces fonctions sont 
intimenent li{\'e}es au monde automorphe).  
Le th{\'e}or{\`e}me principal de cet
article montre que, si $1/q$ est un nombre entier diff{\'e}rent
de $\pm1$ et si $M$ est un nombre impair suffisamment grand, 
alors la dimension
de l'espace vectoriel engendr{\'e} sur $\mathbb Q$ par $1,\zeta_q(3),
\zeta_q(5),\dots, \zeta_q(M)$ est au moins $c_1\cdot\sqrt{M}$, avec
$c_1=0,3358$. Ce r{\'e}sultat peut {\^e}tre consid{\'e}r{\'e} comme un
$q$-analogue du r{\'e}sultat de \cite{ri, br}, 
qui affirme que
la dimension de l'espace vectoriel engendr{\'e} sur $\mathbb Q$ par $1,\zeta(3),
\zeta(5),\dots,\zeta(M)$ est au moins $c_2\cdot\log{M}$, avec
$c_2=0,5906$. Pour les m{\^e}mes valeurs de $q$, une minoration similaire 
pour les valeurs $\zeta_q(s)$ aux entiers $s$ pairs  nous permet de   
red{\'e}montrer un cas particulier d'un r{\'e}sultat de
Bertrand \cite{ber} qui  
affirme la transcendance sur $\mathbb{Q}$ de l'une des deux s{\'e}ries 
d'Eisenstein $E_4(q)$ et $E_6(q)$ pour tout nombre 
complexe $q$ tel que $0<\vert q\vert <1$.

\bigskip

\medskip

\noindent
{\sc Abstract}. We study the arithmetic properties of $q$-analogues
of values $\zeta(s)$ of the Riemann zeta function, 
in particular of the values of the functions
$\zeta_q(s)=
\sum_{k=1} ^{\infty}q^k
\sum_{d\mid k} ^{}d^{s-1}$, $s=1,2,\dots$, where $q$ is a complex number 
with  $\vert q\vert<1$ (these functions are also connected 
with the automorphic world). The main theorem of this article is that,
if $1/q$ is an integer different from $\pm1$, and if $M$
is a sufficiently large odd integer, then the dimension of the vector
space over $\mathbb Q$ which is spanned by $1,\zeta_q(3),
\zeta_q(5),\dots, \zeta_q(M)$ is at least $c_1\cdot\sqrt{M}$,
where $c_1=0,3358$. This result can be regarded as a $q$-analogue of
the result \cite{ri,br} that the dimension of the vector
space over $\mathbb Q$ which is spanned by $1,\zeta(3),
\zeta(5),\dots,\zeta(M)$ is at least $c_2\cdot\log{M}$, with
$c_2=0,5906$. For the same values of $q$, a similar lower bound for the 
values $\zeta_q(s)$ at even integers $s$ provides a new proof of
a special case of a result of Bertrand 
\cite{ber} saying that one of the two Eisenstein series
$E_4(q)$ and $E_6(q)$ is transcendental over $\mathbb Q$ for 
any complex number $q$ such that $0<\vert q\vert <1$.
\end{abstract}

\maketitle

\newpage

\section{Introduction et {\'e}nonc{\'e}s des r{\'e}sultats}
L'{\'e}tude arithm{\'e}tique des valeurs aux entiers $\ge 2$ de la
fonction z{\^e}ta de Riemann est un sujet classique de la th{\'e}orie
analytique des nombres. 
Pour les entiers $2m$ pairs, 
on conna{\^\i}t la formule d'Euler
\begin{equation} \label{eq:Euler} 
\zeta(2m)=(-1)^{m-1}\frac {2^{2m-1}B_{2m}} {(2m)!}\pi^{2m},
\end{equation}
o{\`u} le rationnel $B_m$ est le $m^{\textup{i{\`e}me}}$ nombre de Bernoulli,  
ce qui, joint au fait que $\pi$ est transcendant, {\'e}tablit
que $\zeta(2m)$ l'est aussi pour tout entier $m\ge1$. Concernant les valeurs
$\zeta(2m+1)$, $m\ge1$, notre connaissance est beaucoup plus
faible et pendant longtemps, le seul r{\'e}sultat connu a {\'e}t{\'e} le c{\'e}l{\`e}bre th{\'e}or{\`e}me 
d'Ap{\'e}ry \cite{ap} {\it {\og $\zeta(3)$ est irrationnel\fg}}. R{\'e}cemment, plusieurs 
r{\'e}sultats nouveaux sur les $\zeta(2m+1)$ ont {\'e}t{\'e} d{\'e}montr{\'e}s (voir 
\cite{ri,br,ri1,ri2,zud0} et \cite{fi} pour un survol de ce sujet), 
notamment,  le
th{\'e}or{\`e}me suivant dans \cite{ri, br}.\footnote{Seule est montr{\'e}e 
la minoration aux entiers impairs dans les deux r{\'e}f{\'e}rences
cit{\'e}es~; celle aux entiers pairs est imm{\'e}diate avec la m{\'e}thode utilis{\'e}e.}

\begin{theo} \label{th:rivoal}  
Pour tout entier $A$ pair suffisamment grand, on a
$$
\dim_{\mathbb{Q}} \bigl(\mathbb{Q}+\mathbb{Q}\, 
\zeta(3)+ \mathbb{Q} \,\zeta(5)+ \cdots +
\mathbb{Q}\, \zeta(A-1)\bigr)\ge 
\frac{1+\textup{o}(1)}{1+\log 2}\, \log A
$$
et 
$$
\dim_{\mathbb{Q}} \bigl(\mathbb{Q}+\mathbb{Q}\, 
\zeta(2)+ \mathbb{Q} \,\zeta(4)+ \cdots + 
\mathbb{Q}\, \zeta(A)\bigr)\ge 
\frac{1+\textup{o}(1)}{1+\log 2}\, \log A.
$$
\end{theo}

L'objet de cet article est d'{\'e}tudier les propri{\'e}t{\'e}s
arithm{\'e}tiques de $q$-analogues des valeurs 
de fonction z{\^e}ta de Riemann. En particulier, 
nous consid{\'e}rons les s{\'e}ries $\zeta_q(s)$, 
o{\`u} $s\ge1$ et $q$ est un nombre complexe tels que 
$\vert q\vert <1$, d{\'e}finies par
\begin{equation} \label{eq:zetaq}
\zeta_q(s)=\sum_{k=1}^{\io} \sigma_{s-1}(k)\,q^k=
\sum_{m=1}^{\io} m^ {s-1}\,\frac{q^{m}}{1-q^{m}},
\end{equation}
avec $\sigma_{s-1}(k)=\sum_{d\vert k}d^{s-1}$.
Ces $q$-analogues normalis{\'e}s 
de $\zeta(s)$ ont {\'e}t{\'e} consid{\'e}r{\'e}s dans \cite{jap, zud5}, par
exemple
(voir aussi le
paragraphe~\ref{sec:class} pour la v{\'e}rification que $\zeta_q(s)$
est effectivement un $q$-analogue de $\zeta(s)$). Plus pr{\'e}cis{\'e}ment, notre but est de 
montrer les Th{\'e}or{\`e}mes~\ref{theo:dimietp}, \ref{theo:modu} et 
\ref{theo:357911} suivants. 

\begin{theo} \label{theo:dimietp}
Fixons $q\not =\pm1$ tel que $1/q\in\mathbb{Z}$. Pour tout entier 
$A$ pair, 
on a 
$$
\dim_{\mathbb{Q}} \bigl(\mathbb{Q}+\mathbb{Q}\, 
\zeta_q(3)+ \mathbb{Q} \,\zeta_q(5)+ \cdots +
\mathbb{Q}\, \zeta_q(A-1)\bigr)\ge 
\frac{\pi+\textup{o}(1)}{2\sqrt{\pi^2+12}}\,\sqrt{A}
$$
et 
$$
\dim_{\mathbb{Q}} \bigl(\mathbb{Q}+\mathbb{Q}\, 
\zeta_q(2)+ \mathbb{Q} \,\zeta_q(4)+ \cdots +
\mathbb{Q} \,\zeta_q(A)\bigr)\ge 
\frac{\pi+\textup{o}(1)}{2\sqrt{\pi^2+12}}\,\sqrt{A}.
$$
\end{theo}

La minoration 
de la dimension des valeurs de $\zeta_q(s)$ aux entiers pairs 
poss{\`e}de une int{\'e}\-res\-san\-te traduction en terme de formes modulaires. 
Donnons d'abord 
quelques d{\'e}finitions (voir \cite{jps} pour plus de
d{\'e}tails).  
Posons $q= e^{2i\pi \tau}$, o{\`u} 
$\tau\in\mathcal{H}=\{\tau\in \mathbb{C} : \textup{Im}(\tau)>0\}$.   
Toute matrice 
$\dis\gamma=\binom{a\;\;b}{c\;\;d}\in SL_2(\mathbb{Z})$ agit sur $\mathcal{H}$
par 
$\dis \gamma\cdot\tau=\frac{a\tau+b}{c\tau+d}$. 
Les s{\'e}ries d'Eisenstein $G_{2s}(\tau)$ et $E_{2s}(q)$ sont d{\'e}finies
pour tout 
entier $s\ge 1$ par 
$$
G_{2s}(\tau)=2\zeta(2s)E_{2s}(q)=
\sum_{{(m,n)\in\mathbb{Z}^2\atop (m,n)\neq(0,0)}}\frac{1}{(m\tau+n)^{2s}}
$$
et on a le d{\'e}veloppement en s{\'e}rie de Fourier~:
$$
E_{2s}(q)
=1-\frac{4s}{B_{2s}}\sum_{k=1}^{\io} \sigma_{2s-1}(k)\,q^k.
$$
On a donc 
$E_{2s}(q)=1-\frac{4s}{B_{2s}}\,\zeta_q(2s)$.  
Pour tout $s\ge 2$, les s{\'e}ries d'Eisenstein sont 
modulaires sur le groupe 
$SL_2(\mathbb{Z})$, de poids $2s$, c'est-{\`a}-dire que pour tout $\dis\gamma=
\binom{a\;\;b}{c\;\;d}\in SL_2(\mathbb{Z})$, on a
$G_{2s}(\gamma \cdot\tau)=(c\tau+d)^{2s}G_{2s}(\tau)$.  
La fonction $G_2(\tau)$ n'est pas modulaire mais v{\'e}rifie 
$\dis G_2(\gamma \cdot \tau)=
(a\tau +b)^2G_2(\tau)-\frac{c}{4i\pi}\,\gamma\cdot\tau$. 
On a alors le th{\'e}or{\`e}me suivant.

\begin{theo}
\label{theo:modu}
Pour tout $q\not =\pm1$ tel que $1/q\in\mathbb{Z}$, au moins un 
des deux nombres 
$E_4(q)$ et $E_6(q)$ est transcendant sur $\mathbb{Q}$.
\end{theo}

Ce r{\'e}sultat n'est pas nouveau puisque Bertrand \cite{ber} a montr{\'e}  
que {\it {\og pour tout $q\in\mathbb{C}$ tel que 
$0<\vert q\vert <1$,  au moins un des nombres $E_4(q)$ et $E_6(q)$ est 
transcendant sur $\mathbb{Q}$\fg}}, comme 
cons{\'e}quence  
d'un r{\'e}sultat de Schneider sur les fonctions elliptiques 
\cite[Th{\'e}or{\`e}me~18, p.~64]{sch}. 
Le th{\'e}or{\`e}me st{\'e}phanois \cite{bdgp} pr{\'e}cise le
th{\'e}or{\`e}me de Bertrand~: {\it {\og lorsque $q$ est alg{\'e}brique, 
l'invariant modulaire $J(q)=1728 E_4^3(q)/(E_4^3(q)-E_6^2(q))$ est
transcendant sur $\mathbb Q$\fg}}.
Le r{\'e}sultat d{\'e}finitif dans cette direction 
est le th{\'e}or{\`e}me de Nes\-te\-ren\-ko \cite{ne3}~: 
{\it\og 
pour tout $q\in\mathbb{C}$ tel que 
$0<\vert q\vert <1$, au moins trois des nombres $q$, $E_{2}(q)$, 
$E_4(q)$, $E_6(q)$ sont  alg{\'e}briquement 
ind{\'e}pendants sur $\mathbb{Q}$ \fg}. 
Cependant, la d{\'e}monstration du Th{\'e}or{\`e}me~\ref{theo:modu} est
bas{\'e}e sur une fonction auxiliaire totalement explicite 
(notre $S_n^{[\ep]}(q)$ au paragraphe~\ref{para:resume}), et pas sur
celles, beaucoup 
moins explicites, que l'on peut construire 
avec les outils diophantiens usuels tels le lemme de Siegel ou les
d{\'e}terminants d'interpolation de Laurent 
(voir \cite{bdgp}, \cite{ne3} et \cite{pph} 
pour l'utilisation de ces outils dans le contexte modulaire). 

{\it A contrario}, la minoration de la dimension des valeurs de $\zeta_q(s)$ 
aux entiers impairs est nouvelle 
et ne poss{\`e}de, {\`a} notre connaissance, 
aucune interpr{\'e}tation en terme des formes
modulaires ci-dessus. Cependant, ces valeurs $\zeta_q(s)$, pour $s$ impair, 
ont malgr{\'e} tout un lien avec le monde automorphe, {\it via} les 
s{\'e}ries d'Eisenstein non-holomorphes  
$$
E_s(\tau)=\frac{1}{2}
\sum_{{(m,n)\in\mathbb{Z}^2\atop (m,n)\neq(0,0)}}\frac{y^s}{\vert m\tau+n\vert^{2s}}.
$$ 
Ces s{\'e}ries sont des formes de Maass non paraboliques et 
sont invariantes sous l'action de 
$SL_2(\mathbb{Z})$. Sans rentrer dans les d{\'e}tails, indiquons 
seulement que dans \cite[p. 243]{lewzag}, Lewis et Zagier d{\'e}terminent   
les fonctions 
{\it p{\'e}riodes} $\psi_s^{+}(\tau)$ des $E_s(\tau)$~: pour $s=2\ell+1$, les valeurs $\zeta_q(s)$ sont 
essentiellement {\'e}gales {\`a} $\psi_{\ell+1/2}^{+}(\tau)+\tau^{-2\ell-1}\psi_{\ell+1/2}^{+}(-1/\tau)$. 

Du point de vue diophantien, si l'irrationalit{\'e} de $\zeta_q(1)$ est connue pour diverses 
valeurs de $q$ (voir
\cite{bor} et
les r{\'e}f{\'e}rences donn{\'e}es dans \cite{zud1,zud4}),  celle 
de $\zeta_q(2\ell+1)$  
(pour $\ell\ge 1$) ne l'est pour
aucune valeur de $q$,      
m{\^e}me si la transcendance des ces valeurs 
comme fonctions de $q$ est connue (voir \cite{zud3}). 
Remarquons que pour 
$1/q\in\mathbb Z$, $q\ne\pm1$, des mesures d'irrationalit{\'e} de 
$\zeta_q(1)$ et $\zeta_q(2)$ sont aussi donn{\'e}es par le dernier
auteur dans 
\cite{zud4, zud2} et que Postelmans et Van Assche  \cite{postvan} ont r{\'e}cemment
d{\'e}montr{\'e} l'ind{\'e}pendance lin{\'e}aire de $1$, $\zeta_q(1)$ et $\zeta_q(2)$. 
Enfin, une minoration de dimension dans le contexte diff{\'e}rent 
de la fonction $q$-exponentielle se trouve dans \cite{buv}. 

Compte-tenu de la formule d'Euler \eqref{eq:Euler},
la minoration de la dimension de l'espace des $\zeta$ pairs est  
{\'e}quivalente 
{\`a} la transcendance\footnote{Pour le voir, il suffit d'adapter 
l'argument utilis{\'e} au paragraphe~\ref{sec:demotheomodu}.} de $\pi$~: on peut 
donc consid{\'e}rer que 
les Th{\'e}or{\`e}mes~\ref{theo:dimietp} et \ref{theo:modu} 
sont des $q$-analogues respectifs  du Th{\'e}or{\`e}me~\ref{th:rivoal}
et de la transcendance de $\pi$. 
Comme on le verra au paragraphe~\ref{sec:class}, 
lorsque l'on fait tendre $q$ vers 1, notre construction {\og tend
\fg} vers celle utilis{\'e}e pour montrer le 
Th{\'e}or{\`e}me~\ref{th:rivoal}. 
Puis, on pr{\'e}sente au paragraphe~\ref{sec:zetaq3} des 
$q$-analogues de certaines s{\'e}ries hyperg{\'e}om{\'e}triques 
apparues dans diverses 
d{\'e}monstrations de l'irrationalit{\'e} de $\zeta(3)$. Dans ce contexte, 
notre $q$-analogue de $\zeta(3)$, {\`a} savoir $\zeta_q(3)$, 
appara{\^\i}t naturellement, 
mais, malheureusement, nous ne r{\'e}ussissons {\`a} 
d{\'e}montrer son irrationalit{\'e} pour aucune valeur de~$q$. 

En utilisant les estimations pr{\'e}cises utilis{\'e}es pour montrer le  
Th{\'e}or{\`e}me~\ref{theo:dimietp},  on peut cependant d{\'e}mon\-trer sans effort le th{\'e}or{\`e}me
ci-dessous, de m{\^e}me facture que le r{\'e}sultat du dernier auteur 
\cite{zud0}~: {\it {\og au moins un des nombres $\zeta(5), \zeta(7), \zeta(9), \zeta(11)$ est
irrationnel \fg}}. 

\begin{theo} \label{theo:357911}
Pour tout $q\not =\pm1$ tel que $1/q\in\mathbb{Z}$, au moins un des
nombres $\zeta_q(3),\zeta_q(5),\break 
\zeta_q(7),\zeta_q(9),
\zeta_q(11)$ est irrationnel.
\end{theo}

Il est probable que l'on puisse am{\'e}liorer ce r{\'e}sultat 
(i.e., d{\'e}montrer qu'au moins un des nombres $\zeta_q(3),\zeta_q(5),
\zeta_q(7),\zeta_q(9)$ est irrationnel) en 
utilisant une s{\'e}rie l{\'e}g{\`e}rement diff{\'e}rente de notre s{\'e}rie
$S_n^{[\ep]}(q)$ ci-dessous (voir la Remarque~(2) au 
paragraphe~\ref{para:resume}), {\`a} condition que l'on puisse d{\'e}montrer une certaine {\og 
conjecture des d{\'e}nominateurs\fg} pour cette s{\'e}rie (voir la remarque 
{\`a} la fin de l'article).

\section{D{\'e}monstration du Th{\'e}or{\`e}me \ref{theo:modu}}
\label{sec:demotheomodu}
Le Th{\'e}or{\`e}me~\ref{theo:modu} est 
une cons{\'e}quence du r{\'e}sultat suivant, dont on trouvera la 
d{\'e}mon\-stra\-tion dans \cite{jps}.

\begin{prop}[\sc Structure de l'espace des formes modulaires] Soit $M(q)$ une
forme modulaire sur $SL_2(\mathbb{Z})$, holomorphe sur le disque   
$\{q\in \mathbb{C} : \vert q\vert  <1\}$,  de poids $n\ge 1$. Alors $n$ est pair et il 
existe des nombres complexes 
$c_{a,b}$ tels que 
$
\dis M(q)=\sum c_{a,b}E_4(q)^aE_6(q)^b
$,   
o{\`u} la sommation est {\'e}tendue {\`a} tous les couples d'entiers 
positifs $(a,b)$ tels que $4a+6b=n$.
\end{prop}

Lorsque $M(q)$ est une s{\'e}rie d'Eisenstein $E_{2n}(q)$, les  $c_{a,b}$ 
peuvent {\^e}tre pris rationnels. Comme pour tout entier pair $j\ge 4$, 
les formes modulaires $\zeta_q(j)$ de poids $j$ v{\'e}rifient 
les hypoth{\`e}ses de 
ce th{\'e}or{\`e}me, on en d{\'e}duit  
que pour tout $A$ pair, le $\mathbb{Q}$-espace vectoriel $V(A)$ engendr{\'e} 
par $1, \zeta_q(4), 
\zeta_q(6), \ldots, \zeta_q(A)$ est inclus dans le $\mathbb{Q}$-espace 
vectoriel $W(A)$ 
engendr{\'e} par les puissances 
$E_4(q)^aE_6(q)^b$, avec $0\le a\le A/4$ et $0\le b\le A/6$. 
Fixons $q\not =\pm1$ tel que $1/q\in\mathbb{Z}$~; si les nombres $E_4(q)$ 
et $E_6(q)$ {\'e}taient 
tous les deux alg{\'e}briques sur $\mathbb{Q}$, la dimension de $W(A)$, et donc 
celle de $V(A)$, serait born{\'e}e 
ind{\'e}pendamment de $A$. Or ceci est impossible puisque, pour $A$ assez grand, 
cela contredit 
la minoration de la dimension de $V(A)$ fournit par le
Th{\'e}or{\`e}me~\ref{theo:dimietp}. D'o{\`u} le 
Th{\'e}or{\`e}me~\ref{theo:modu}.

\begin{remark} 
D'apr{\`e}s le th{\'e}or{\`e}me de Nesterenko \cite{ne3}, les nombres 
$\zeta_q(2)$, $\zeta_q(4)$ et 
$\zeta_q(6)$ sont alg{\'e}bri\-que\-ment ind{\'e}pendants pour 
$q\not =\pm1$ tel que $1/q\in\mathbb{Z}$. 
Pour tout $A\ge 2$ pair, on a donc en fait
$$
\dim_{\mathbb{Q}} \bigl(\mathbb{Q}+\mathbb{Q}\,
\zeta_q(2)+ 
\mathbb{Q} \,\zeta_q(4)+ \cdots +
\mathbb{Q} \,\zeta_q(A)\bigr)= A/2+1.
$$
\end{remark}

\section{D{\'e}monstrations des Th{\'e}or{\`e}mes \ref{theo:dimietp} et 
\ref{theo:357911}}
\label{sec:demodimietp}
\subsection{La m{\'e}thode}\label{para:resume}
La d{\'e}monstration est bas{\'e}e sur la construction de combinaisons 
li\-n{\'e}\-aires rationnelles 
en 1 et les $\zeta_q(j)$, avec $j$ pair ou impair, exclusivement. 
Pour cela, introduisons  
la s{\'e}rie  
$$
S_n(q)
= (q;q)_n^{A-2r}
\sum_{k=1}^{\io}q^k 
\frac{(q^{k-rn};q)_{rn}\,(q^{k+n+1};q)_{rn}}{(q^k;q)_{n+1}^{A}}\;q^{(k-1/2)(A-2r)n/2},
$$
o{\`u} les $q$-factoriels
montants sont d{\'e}finis par $(\alpha;q)_m=(1-\alpha)(1-q\alpha)\cdots
(1-q^{m-1}\alpha)$ si $m\ge1$, et $(\alpha;q)_0=1$, et
o{\`u} l'on suppose $A$, $n$, $r$ entiers tels que 
\begin{equation}
n\ge 0,\; A \;\textup{pair et}\;  1\le r\le  A/2.
\label{eq:condanr}
\end{equation}
Cette s{\'e}rie, qui   
converge pour tout nombre complexe $q$ tel que $\vert q\vert \ne 1$, est 
une s{\'e}rie hyperg{\'e}om{\'e}trique 
basique (voir le paragraphe~\ref{ssec:combE} pour plus de d{\'e}tails).   
Fixons $\ep\in\{0,1\}$ et $A$, $n$, $r$ v{\'e}rifiant~\eqref{eq:condanr}, 
et $q$ un nombre complexe tel 
que $\vert q\vert <1$. La  s{\'e}rie $S_n(q)$ v{\'e}rifie les propri{\'e}t{\'e}s
suivantes~:

\medskip
\begin{itemize}
\item[A)] Il existe des fractions rationnelles $P_{s,n}^{[\ep]}(q)$ de $\mathbb{Q}(q)$, 
avec $s\in\{2,\ldots, A\}$,  
telles que 
$$
S_n^{[\ep]}(q)=S_n(q)+(-1)^{\ep}q^{-n}S_n(1/q)=
P_{0,n}^{[\ep]}(q)+
\underset{s\equiv \ep \,(\textup {mod }2)}{\sum_{{s=2}}^{A}} 
P_{s,n}^{[\ep]}(q)\,\zeta_q(s).
$$

\item[B)] On a 
$\dis 
\lim_{n\to +\io}\frac{1}{n^2}\log\left\vert S_{n}^{[\ep]}(q)\right\vert =- \frac 12
r(A-2r) \log\vert 1/q\vert \ne 0$.

\item[C)] Pour tout $s\in\{0,\ldots, A\}$,  on a 
$\dis 
\limsup_{n\to +\io}\frac{1}{n^2}\log\vert P_{s,n}^{[\ep]}(q)\vert \le
\frac {1} {8}\left( A+4r^2\right) \log \vert 1/q\vert
$.

\item[D)] Il existe $D_n(q)\in\mathbb{Q}(q)$ tel que pour $s\in\{0,\ldots, A\}$ on a 
$D_n(q)P_{s,n}^{[\ep]}(q)\in\mathbb{Z}[1/q]$ et 
\begin{equation} \label{eq:Dnq}
\lim_{n\to +\io} \frac{1}{n^2}\log \vert D_n(q) \vert
= \left(\frac{3}{\pi^2}A+\frac A8 + \frac {r^2} 2\right)\log \vert 1/q\vert.
\end{equation}
\end{itemize}

\begin{remarks} (1) La d{\'e}monstration du point~D)
montrera que l'estimation arithm{\'e}tique 
donn{\'e}e dans ce point est optimale, i.e.\ que l'on
ne peut pas remplacer ce $D_n(q)$ par un  $\tilde
D_n(q)$ ayant les m{\^e}mes propri{\'e}t{\'e}s que 
$D_n(q)$, mais tel que la limite de  
$\frac{1}{n^2}\log \vert \tilde D_n(q) \vert$ soit
plus petite que le membre de droite de \eqref{eq:Dnq}. 
Ceci s'applique 
probablement aussi {\`a} l'estimation donn{\'e}e au point~C) (voir la 
remarque {\`a} la fin du paragraphe~\ref{sec:P}). 
Par cons{\'e}quent, la minoration des dimensions donn{\'e}e au
Th{\'e}or{\`e}me~\ref{theo:dimietp} est tr{\`e}s certainement la meilleure 
que l'on puisse obtenir par {\it la m{\'e}thode 
utilis{\'e}e} et {\it notre s{\'e}rie $S_{n}^{[\ep]}(q)$}. 

\medskip

(2) En revanche, il est sans doute possible d'am{\'e}liorer 
certains des points A) {\`a} D) en
utilisant {\it d'autres s{\'e}ries} et en {\it raffinant la m{\'e}thode}. 
Par exemple, on peut obtenir les m{\^e}mes r{\'e}sultats A)--D) en
utilisant la s{\'e}rie
\begin{equation} \label{eq:Salt}
\tilde S_n(q)
= (q;q)_n^{A-2r}
\sum_{k=1}^{\io} \left(1-q^{2k+n}\right)
\frac{(q^{k-rn};q)_{rn}\,(q^{k+n+1};q)_{rn}}
{(q^k;q)_{n+1}^{A}}\;q^{k((A-2r)n/2+A/2-1)}
\end{equation}
en lieu et place de $S_{n}^{[\ep]}(q)$ et  en appliquant des arguments 
similaires. Or cette s{\'e}rie a semble-t-il certains avantages arithm{\'e}tiques
sur la s{\'e}rie $S_{n}^{[\ep]}(q)$ (voir la remarque {\`a} la 
fin de l'article).
\end{remarks}

\medskip
Le Th{\'e}or{\`e}me~\ref{theo:dimietp} r{\'e}sulte des
propri{\'e}t{\'e}s A)--D) ci-dessus et du cas particulier  
suivant d'un crit{\`e}re d'in\-d{\'e}\-pen\-dance lin{\'e}aire, d{\^u}
{\`a} Nesterenko \cite{ne1}.  

\medskip
\begin{prop}[\sc Crit{\`e}re d'ind{\'e}pendance lin{\'e}aire]
\label{prop:Nest}
Soient un 
entier $N\ge 2$ et $\vartheta_1$, \dots,
$\vartheta_N$ des r{\'e}els. Supposons qu'il existe $N$ suites d'entiers
$(p_{j,n})_{n\ge 0}$ et des r{\'e}els $\al$ et $\be$ avec $\beta>0$ tels que 

\begin{itemize}
\item[\em i)] 
pour tout $j\in\{1, \ldots, N\}$, on a $\dis \limsup_{n\to
+\io}\frac{1}{n^2} \log \vert
p_{j,n}\vert \le \be$~;
\item[\em ii)]
$\dis \lim_{n\to
+\io}\frac{1}{n^2} \log \vert
p_{1,n}\vartheta_1+\cdots +p_{N,n}\vartheta_N\vert = \al$.
\end{itemize}

\medskip
\noindent Alors la dimension  du $\mathbb{Q}$-espace vectoriel
engendr{\'e} par $\vartheta_1$, \dots, $\vartheta_N$ est plus grande que 
$$
1-\frac{\al}{\be}.
$$  
\end{prop}

\noindent {\it D{\'e}monstration du Th{\'e}or{\`e}me \ref{theo:dimietp}.} 
Supposons $q\not =\pm1$ tel que $1/q\in\mathbb{Z}$ et posons
$$
{\mathfrak
P}_{s,n}^{[\ep]}(q)=D_n(q)P_{s,n}^{[\ep]}(q)\in\mathbb{Z}[1/q]\subset \mathbb{Z}
\quad \textup{et}\quad {\mathfrak S}_n^{[\ep]}(q)=D_n(q)S_{n}^{[\ep]}(q).
$$
En vertu du point~A) ci-dessus, on a
$$
{\mathfrak S}_n^{[\ep]}(q)=
{\mathfrak P}_{0,n}^{[\ep]}(q)+
\underset{s\equiv \ep\,(\textup {mod }2)}{\sum_{{s=2}}^{A}}
{\mathfrak P}_{s,n}^{[\ep]}(q)\,\zeta_q(s).
$$
Gr{\^a}ce aux estimations r{\'e}sum{\'e}es par les 
point~B) {\`a} D) ci-dessus, 
la Proposition~\ref{prop:Nest} 
appliqu{\'e}e {\`a} la combinaison 
lin{\'e}aire ${\mathfrak S}_n^{[\ep]}(q)$ 
permet alors d'affirmer que pour tout $A\ge 2$ pair, les deux dimensions 
$$
\dim_{\mathbb{Q}} \bigl(\mathbb{Q}+\mathbb{Q}\,
\zeta_q(3)+ 
\mathbb{Q} \,\zeta_q(5)+ \cdots +
\mathbb{Q}\, \zeta_q(A-1)\bigr)
$$
et 
$$
\dim_{\mathbb{Q}} \bigl(\mathbb{Q}+\mathbb{Q}\,
\zeta_q(2)+ 
\mathbb{Q} \,\zeta_q(4)+ \cdots +
\mathbb{Q} \,\zeta_q(A)\bigr)
$$
sont  minor{\'e}es, pour tout $r\in \{1,\ldots, A/2 \}$, par 
\begin{equation} \label{eq:delta} 
\delta(A,r)=\frac{4rA+A-4r^2}{\left(\frac{24}{\pi^2}+2\right)A+8r^2}.
\end{equation}
On choisit alors $r$ de la forme $u\sqrt{A}$ o{\`u} $u>0$~; un calcul
imm{\'e}diat montre que 
$$
\max_{u>0}\lim_{A\to+\io}\frac{1}{\sqrt{A}}\;\delta(A,u\sqrt{A})=
\max_{u>0}\left(\frac{4u}{\frac{24}{\pi^2}+2+8u^2}\right)=
\frac{\pi}{2\sqrt{\pi^2+12}}\approx 
0,335892.
$$
Le Th{\'e}or{\`e}me~\ref{theo:dimietp} en d{\'e}coule. $\square$

\medskip
\noindent {\it D{\'e}monstration du Th{\'e}or{\`e}me \ref{theo:357911}.} 
On r{\'e}p{\`e}te la d{\'e}monstration du Th{\'e}or{\`e}me~\ref{theo:dimietp}.
En choisissant $A=12$ et $r=2$ dans~\eqref{eq:delta}, 
on obtient $\delta(12,2)=1,080059\ldots>1$. La dimension
de l'espace vectoriel engendr{\'e} sur $\mathbb{Q}$ par $1,\zeta_q(3),\zeta_q(5),
\zeta_q(7),\zeta_q(9),
\zeta_q(11)$ est donc au moins 2. $\square$ 

\medskip
Les points A), B), C) et D)  sont d{\'e}montr{\'e}s respectivement aux 
Lemmes~\ref{lem:dichoto},~\ref{lem:asympserie},~\ref{lem:asympoly} et~\ref{lem:arith}, dans 
les sous-paragraphes qui suivent. On fait la convention 
suivante~: la valeur d'une fonction
$F(z;q)$ en $z=1$ sera not{\'e}e indiff{\'e}remment  $F(1;q)$ ou $F(q)$,
s'il n'y a pas d'ambigu{\"\i}t{\'e} possible.

\subsection{Des fonctions interm{\'e}diaires} 
On suppose dans tout ce paragraphe que $\vert q\vert <1$. Pour tout $s \ge 1$, la fonction 
$\mathcal{Z}_s(z;q)$ est d{\'e}finie par la s{\'e}rie 
$$
\mathcal{Z}_{s}(z;q)=\sum_{k=1}^{\io}\frac{q^k}{(1-q^k)^s}\,z^{-k},
$$ 
qui converge pour tout $z$ tel que 
$\vert q\vert <\vert z\vert $ et en particulier en $z=1$. 
On aura besoin de consid{\'e}rer 
{\it simultan{\'e}ment} les fonctions $\mathcal{Z}_s(z;q)$ et 
$\mathcal{Z}_s(1/z;1/q)$ : si $s\ge 2$, alors 
cette derni{\`e}re s{\'e}rie est convergente d{\`e}s que $\vert z\vert <\vert q\vert^{1-s}$ 
et il est donc possible de d{\'e}finir  simultan{\'e}ment nos deux  fonctions 
sur $\vert q\vert <\vert z\vert <\vert q\vert^{1-s}$ et en particulier toujours en $z=1$. 
Lorsque $s=1$, la s{\'e}rie  $\mathcal{Z}_1(1/z;1/q)$  converge sur 
$\vert z\vert <\vert q\vert^{1-s}=1$, 
et les deux fonctions $\mathcal{Z}_1(z;q)$ et $\mathcal{Z}_1(1/z;1/q)$ sont 
donc d{\'e}finies 
simultan{\'e}ment sur $\vert q\vert <\vert z\vert < 1$. 

\medskip
Sauf mention contraire, on supposera  dans toute la suite de cet article que 
$\vert q\vert <\vert z\vert < 1$. On fera  
tendre $z$ vers $1$ pour  d{\'e}finir une expression en $z=1$. 
\medskip

Avec cette condition, toutes les s{\'e}ries de ce paragraphe sont 
absolument convergentes et 
en d{\'e}veloppant  $1/(1- q^k)^s$ en
s{\'e}rie de puissances de $q^k$,  on obtient 
pour tout 
$s\ge 1$ que
\begin{equation} \label{eq:Zs}
\mathcal{Z}_s(z;q)=
\frac{1}{(s-1)!}\sum_{k= 1}^{\io}
\sum_{\ell=1}^{\io}\ell(\ell+1)\cdots (\ell+s-2)z^{-k} q^{k\ell} 
\end{equation}
et  
\begin{equation} \label{eq:Zs1/q}
\mathcal{Z}_s(1/z;1/q)=
\frac{(-1)^s}{(s-1)!}\sum_{k= 1}^{\io}
\sum_{\ell=0}^{\io}
\ell(\ell-1)\cdots (\ell-s+2)z^k q^{k\ell},
\end{equation}
en convenant que pour $s=1$ les produits vides
$\ell(\ell+1)\cdots(\ell+s-2)$ et 
$\ell(\ell-1)\cdots(\ell-s+2)$ sont interpr{\'e}t{\'e}s comme 1.
Les fonctions $\mathcal{Z}_{s}(z;q)$ sont donc des combinaisons 
lin{\'e}aires
en les fonctions 
$\mathcal{L}_j(z;q)$, 
d{\'e}finies par les s{\'e}ries (avec $j\ge1$)~:
\begin{equation} \label{eq:Lq}
\mathcal{L}_j(z;q)=\sum_{k=1}^{\io}\sum_{\ell=1}^{\io}\ell^{j-1} z^{-k} q^{k\ell} = 
\sum_{k=1}^{\io}\sigma_{j-1}(z;k)\,q^k \quad
\textup{avec} \quad \sigma_{j-1}(z;k)=\sum_{d\vert k}d^{j-1}z^{-k/d}.
\end{equation}
Plus pr{\'e}cis{\'e}ment, si l'on utilise les nombres de Stirling de premi{\`e}re
esp{\`e}ce (sans signe) $c(N,j)$, d{\'e}finis par (voir 
\cite[Sec.~1.3]{StanAP})
\begin{equation} \label{eq:Stirling1}
\ell(\ell+1)\cdots(\ell+N-1)=
\sum _{j=1} ^{N}c(N,j)\ell^j,
\end{equation}
on a
\begin{equation} \label{eq:Z->E}
\mathcal Z_s(z;q)=\frac {1} {(s-1)!}
\sum _{j=2} ^{s}c(s-1,j-1)\,\mathcal L_j(z;q).
\end{equation}
{\'E}videmment, le coefficient de $\mathcal L_j(z;q)$ est un nombre
rationnel ne d{\'e}pendant que de $s$ et $j$.
De plus, en comparant les d{\'e}finitions \eqref{eq:zetaq} et
\eqref{eq:Lq}, on remarque que $\mathcal 
L_j(1;q)=\zeta_q(j)$.

\begin{lem}\label{lem:crucial}
{\em i)} Pour tout entier $s\ge 2$, on a 
\begin{equation} \label{eq:Z-Z->E} 
\mathcal{Z}_s(q)-\mathcal{Z}_s(1/q)=\frac {2} {(s-1)!}
\underset{j\textup{ impair}}{\sum _{j=3} ^{s}}c(s-1,j-1)\,\zeta_q(j).
\end{equation}

{\em ii)} Pour tout entier $s\ge 2$, on a 
\begin{equation} \label{eq:Z+Z->E}
\mathcal{Z}_s(q)+\mathcal{Z}_s(1/q)=\frac {2} {(s-1)!}
\underset{j\textup{ pair}}{\sum _{j=2} ^{s}}c(s-1,j-1)\,\zeta_q(j).
\end{equation}

{\em iii)} On a 
$\displaystyle \mathcal{Z}_{1}(z;q)=\mathcal{L}_{1}(z;q)$ et 
$\dis \mathcal{Z}_{1}(1/z;1/q)=-\mathcal{L}_{1}(1/z;q)-\frac{z}{1-z}$.
\end{lem}

\begin{proof}
i) et ii) Supposons $s\ge 2$. 
Comme les s{\'e}ries  $\mathcal{Z}_s(z;q)$ et $\mathcal{Z}_s(1/z;1/q)$ 
sont convergentes en $z=1$, on a, en vertu de \eqref{eq:Zs} et \eqref{eq:Zs1/q},
$$
\mathcal{Z}_s(q)+(-1)^{\ep}\mathcal{Z}_s(1/q)=\frac{1}{(s-1)!} 
\sum_{k= 1}^{\io}\sum_{\ell=1}^{\io}\bigl((\ell)_{s-1}+
(-1)^{1+\ep}
(-\ell)_{s-1}\bigr)q^{k\ell},
$$
o{\`u} le symbole de Pochhammer $(\alpha)_m$ est d{\'e}fini par
$(\alpha)_m=\alpha(\alpha+1)\cdots (\alpha+m-1)$, $m\ge1$, et $(\alpha)_0=1$. 
Il est justifi{\'e} de d{\'e}buter la sommation sur $\ell$ {\`a} partir de 1
puisque pour $s\ge2$, on a $(0)_{s-1}=0$.
Si $\ep=0$,
respectivement 1, alors 
$(\ell)_{s-1}+(-1)^{1+\ep}
(-\ell)_{s-1}$ est un polyn{\^o}me impair, respectivement pair. Cela
d{\'e}montre i) et ii), sauf qu'il reste l'{\'e}ventualit{\'e} que 
$\zeta_q(1)$ apparaisse aussi dans la combinaison
lin{\'e}aire pour $\ep=1$~: comme $(0)_{s-1}=0$, c'est en fait impossible. 
Les formes explicites \eqref{eq:Z-Z->E} et \eqref{eq:Z+Z->E} 
r{\'e}sultent de \eqref{eq:Stirling1}.

\medskip
iii) En 
utilisant \eqref{eq:Zs}, on obtient
$$
\mathcal{Z}_{1}(z;q)=
\sum_{k= 1}^{\io}\sum_{\ell=1}^{\io}z^{-k} q^{k\ell}=\mathcal{L}_{1}(z;q).
$$
De m{\^e}me, en utilisant \eqref{eq:Zs1/q}, on obtient
$$
\mathcal{Z}_{1}(1/z;1/q)=-\sum_{k= 1}^{\io}
\sum_{\ell=0}^{\io}z^{k} q^{k\ell}=
-\sum_{k=1}^{\io}\sum_{\ell=1}^{\io}z^{k} q^{k\ell}-
\sum_{k=1}^{\io}z^k=-\mathcal{L}_{1}(1/z;q)-\frac{z}{1-z}.
$$
\end{proof}

\subsection{Construction de combinaisons lin{\'e}aires 
en les $\zeta_q(s)$ }\label{ssec:combE}
Introduisons la fraction rationnelle en $T$~: 
\begin{equation*}
R_n(T;q)
= (q;q)_n^{A-2r}q^{-(A-2r)n/4}\,
\frac{(q^{-rn}T;q)_{rn}\,(q^{n+1}T;q)_{rn}}{(T;q)_{n+1}^{A}}\;T^{(A-2r)n/2}
\end{equation*}
avec $n\ge 0$, $A$ pair et $1\le r\le  A/2$.
D{\'e}finissons {\'e}galement la
s{\'e}rie~:
$$
S_n(z;q)=\sum_{k=1}^{\io}q^k R_n(q^k;q)z^{-k},
$$
qui converge (au moins) pour tout $z$ et tout $q$ tels que 
$\vert q\vert <\vert z\vert \le 1$ 
(c'est ici que l'on utilise la condition $r\le A/2$)~: 
en $z=1$, cette s{\'e}rie co{\"\i}ncide avec 
la s{\'e}rie $S_n(q)$ introduite au paragraphe~\ref{para:resume}. 
En utilisant la notation classique
\begin{equation} \label{eq:hyp}
{}_r\phi_s\!\left[\begin{matrix} a_1,\dots,a_r\\
b_1,\dots,b_s\end{matrix}; q,
z\right]=\sum _{n=0} ^{\infty}\frac {\poq{a_1}{n}\cdots\poq{a_r}{n}}
{\poq{q}{n}\poq{b_1}{n}\cdots\poq{b_s}{n}}\left((-1)^nq^{\binom
n2}\right)^{s-r+1}z^n
\end{equation}
pour la s{\'e}rie hyperg{\'e}om{\'e}trique basique (voir \cite{gr}), on peut {\'e}crire
$S_n(z;q)$ comme
\begin{multline*}
q^{1 + \frac{A n}{4} + \frac{rn }{2} + \frac{rA n^2 }{2} -
r^2 n^2 }
\frac{
({\let \over / \def\frac#1#2{#1 / #2} q}; q) _{n} ^{A - 2 r}\,
({\let \over / \def\frac#1#2{#1 / #2} q}; q) _{r n} \,
({\let \over / \def\frac#1#2{#1 / #2} q^ {(r+1)n+2}}; q)
_{r n}  }{{z^{ rn+1}(\let \over /
\def\frac#1#2{#1 / #2} q^ { rn+1}}; q) _{n+1} ^A}\\
\times
{}_{A+2}\phi_{A+1}\!\left[\begin{matrix} 
q^ { (2r+1)n+2}, q^ { r n+1},\dots,q^ {r n+1}\\
q^ {(r+1)n+2},\dots,q^ {(r+1)n+2}\end{matrix}; q, 
\frac {q^ { \frac{1}{2}(A - 2 r)n+1}} z\right].
\end{multline*}

D{\'e}composons 
$R_n(T;q)$ en {\'e}l{\'e}ments simples~: 
\begin{align} \label{eq:wadim}
R_n(T ; q)
&=q^{-An(n+1)/2}(q;q)_n^{A-2r}q^{-(A-2r)n/4}\\
\notag
&\quad
\times \frac{(1-q^{-rn}T)\cdots(1-q^{-1}T)\cdot(1-q^{n+1}T)\cdots
(1-q^{(r+1)n}T)}{(T-1)^{A}\cdots (T-q^{-n})^{A}}\,T^{(A-2r)n/2}\\
\notag
&= \sum_{s=1}^A\sum_{j=0}^n \frac{c_{s,j,n}(q)}{(T-q^{-j})^{s}}=
\sum_{s=1}^A\sum_{j=0}^n \frac{d_{s,j,n}(q)}{(1-q^jT)^{s}}
\end{align}
avec $d_{s,j,n}(q)=(-1)^sq^{js}c_{s,j,n}(q)$. Il n'y a pas de 
partie principale car la diff{\'e}rence des
degr{\'e}s du num{\'e}rateur et du d{\'e}nominateur de $R_n(T;q)$
est $-A-(A-2r)n/2$,
qui est strictement n{\'e}gative puisque 
$1\le r\le A/2$~; de plus,  selon la formule usuelle, on a
\begin{equation} \label{eq:c} 
c_{s,j,n}(q)=\frac{1}{(A-s)!}
\frac{\dd^{A-s}}{\dd T^{A-s}}\bigl(R_n(T;q)(T-q^{-j})^A\bigr)
\big\vert_{T=q^{-j}}.
\end{equation}
On construit alors des combinaisons lin{\'e}aires en
les fonctions $\mathcal{Z}_s(z;q)$ de la fa{\c c}on suivante~:
\begin{eqnarray*}
S_n(z;q)&=&\sum_{s=1}^A\sum_{j=0}^n
d_{s,j,n}(q)\,q^{-j}z^j\sum_{k=1}^{\io}\frac{q^{k+j}}{(1-q^{k+j})^s}z^{-k-j}\\
&=&P_{0,n}(z;q)+\sum_{s=1}^A P_{s,n}(z;q)\,\mathcal{Z}_s(z;q)
\end{eqnarray*}
avec, pour $s\ge 1$, 
$$
P_{s,n}(z;q)=\sum_{j=0}^n d_{s,j,n}(q)\,q^{-j}z^j
\quad \textup{et}\quad
P_{0,n}(z;q)=-\sum_{s=1}^A
\sum_{j=1}^n
\sum_{k=1}^j d_{s,j,n}(q)\,q^{k-j}
\frac{z^{j-k}}{(1-q^{k})^s}.
$$
Jusqu'{\`a} pr{\'e}sent, nous n'avons pas exploit{\'e} la forme particuli{\`e}re
du num{\'e}rateur de $R_n(T;q)$, ce que nous allons maintenant faire. 

\begin{lem}\label{lem:sympoly}
Soit $A$ pair.
Pour tout $s\in\{ 1, \ldots, A\}$, on a la relation de
r{\'e}ciprocit{\'e}~: 
$$
z^nq^{-n}P_{s,n}(1/z;1/q)=P_{s,n}(z;q).
$$
\end{lem}
\begin{proof} On v{\'e}rifie imm{\'e}diatement que 
$$
R_n(q^nT;1/q)=R_n(T;q). 
$$
Par unicit{\'e} du d{\'e}veloppement en {\'e}l{\'e}ments simples
\eqref{eq:wadim}, cette sym{\'e}trie de 
$R_n(T ; q)$ se traduit par 
$$
d_{s,n-j,n}(q)=d_{s,j,n}(1/q),
$$
ce qui prouve le lemme.
\end{proof}

En utilisant l'identit{\'e} triviale 
$(1/q^{\alpha};1/q)_{\beta}=(-1)^{\beta}q^{-\al\be -
\be(\be-1)/2}(q^{\al};q)_{\be}$, on montre sans difficult{\'e} que 
$$
S_n(1/z;1/q)=q^{An/2}S_n(q^{2-A}/z;q),
$$
qui est convergente sur  
$\vert q\vert <\vert z\vert \le 1$~: on peut donc consid{\'e}rer 
simultan{\'e}ment les s{\'e}ries 
$S_n(z;q)$ et $S_n(1/z;1/q)$. Cela va nous permettre de 
construire une nouvelle combinaison lin{\'e}aire r{\'e}alisant la 
dichotomie attendue entre les valeurs de la fonction $\zeta_q(s)$ aux entiers $s$ pairs et impairs. 

\begin{lem}\label{lem:dichoto}
Soit $A$ pair et $\ep\in\{0,1\}$. 
Il existe des fractions rationnelles $P_{s,n}^{[\ep]}(q)$ de $\mathbb{Q}(q)$, 
avec $s\in\{0,\ldots, A\}$,  
telles que 
$$
S_n^{[\ep]}(q)=S_n(q)+(-1)^{\ep}q^{-n}S_n(1/q)=
P_{0,n}^{[\ep]}(q)+
\underset{s\equiv \ep\,(\textup {mod }2)}{\sum_{{s=2}}^{A}} 
P_{s,n}^{[\ep]}(q)\,\zeta_q(s).
$$
\end{lem}
On explicitera les  $P_{s,n}^{[\ep]}(q)$ au cours de la d{\'e}monstration. 

\begin{proof} On se concentre sur le cas $\ep=1$, le cas  $\ep=0$ {\'e}tant
  similaire. Compte-tenu de la relation de r{\'e}ciprocit{\'e} mise en {\'e}vidence au
Lemme~\ref{lem:sympoly}, on a 
\begin{multline*}
S_n(z;q)-q^{-n}z^nS_n(1/z;1/q)
\\=
P_{0,n}(z;q)-q^{-n}z^nP_{0,n}(1/z;1/q)+
\sum_{s=1}^A P_{s,n}(z;q)\bigl(
\mathcal{Z}_s(z;q)-\mathcal{Z}_s(1/z;1/q)\bigr). 
\end{multline*}
Les s{\'e}ries 
$S_n(z;q)$ et $q^nz^nS_n(1/z;1/q)$ sont convergentes en $z=1$, ainsi que les fonctions 
$\mathcal{Z}_s(z;q)-\mathcal{Z}_s(1/z;1/q)$ pour $s\ge 2$. Le seul terme 
potentiellement divergent est 
$\mathcal{Z}_1(z;q)-\mathcal{Z}_s(1/z;1/q)$~:  pour contrer cette divergence, 
le polyn{\^o}me (en la variable $z$) 
$P_{1,n}(z;q)$ s'annule  n{\'e}cessairement en $z=1$. 
Or en utilisant le point~iii) du Lemme~\ref{lem:crucial}, on a 
\begin{multline*}
\lim_{{z\to 1^-}} P_{1,n}(z;q)\bigl(
\mathcal{Z}_1(z;q)-\mathcal{Z}_1(1/z;1/q)\bigr)\\=
\lim_{{z\to 1^-}} P_{1,n}(z;q)\bigl(
\mathcal{L}_1(z;q)+\mathcal{L}_1(1/z;q)\bigr)+
\lim_{{z\to 1^-}}
\frac{z\,P_{1,n}(z;q)}{1-z}=-\frac{\partial P_{1,n} }{\partial z}(1;q)
\end{multline*}
puisque $\mathcal{L}_1(1;q)$ est convergente. 
Par ailleurs, le point~i) du m{\^e}me Lemme~\ref{lem:crucial} montre que
pour $s\ge 2$, 
$$
 \mathcal{Z}_s(q)-\mathcal{Z}_s(1/q)=
\underset{j\,\textup{impair}}{\sum_{{j=3}}^{s}}
\alpha_{s,j}\, \zeta_q(j)
$$
o{\`u} les coefficients $\al_{s,j}$ sont des rationnels ind{\'e}pendants de
$q$ et dont un d{\'e}nominateur commun est 
$(s-1)!$ 
(ce fait nous servira au
cours de la d{\'e}monstration du Lemme~\ref{lem:arith}~; 
en fait,
$\alpha_{s,j}=2c(s-1,j-1)/(s-1)!$, o{\`u} $c(s-1,j-1)$ est un nombre de
Stirling de premi{\`e}re esp{\`e}ce sans signe). 
Donc 
$$
S_n^{[1]}(q)=P_{0,n}(1;q)-q^{-n}P_{0,n}(1;1/q)-
\frac{\partial P_{1,n} }{\partial z}(1;q)+
\sum_{s=3}^{A}\underset{j\,\textup{impair}}
{\sum_{{j=3}}^{s}}\alpha_{s,j}P_{s,n}(1;q)\,\zeta_q(j).
$$
Il suffit donc de poser
$
\displaystyle P_{0,n}^{[1]}(q)=P_{0,n}(1;q)-q^{-n}P_{0,n}(1;1/q)-
\frac{\partial P_{1,n} }{\partial z}(1;q)
$
et, pour $s$ impair de $\{3, \ldots, A-1$\}, 
\begin{equation} \label{eq:S->P} 
\displaystyle P_{s,n}^{[1]}(q)=
\sum_{k=s}^{A}\alpha_{k,s}P_{k,n}(1;q).
\end{equation}

Dans le cas $\ep=0$, on trouve  
$\displaystyle P_{0,n}^{[0]}(q)=P_{0,n}(1;q)+q^{-n}P_{0,n}(1;1/q)+
\frac{\partial P_{1,n} }{\partial z}(1;q)$ 
et, pour $s$ pair de $\{2, \ldots, A$\}, 
$\displaystyle P_{s,n}^{[0]}(q)=
\sum_{k=s}^A\alpha_{k,s}P_{k,n}(1;q).$
\end{proof}

Lorsque l'on essaye de relier la s{\'e}rie $S^{[\ep]}_n(q)$ {\`a} la notation
classique \eqref{eq:hyp}, on obtient
\begin{align}
\label{eq:Sexplicit}
S^{[\ep]}_n&(q)=
\sum _{k=1}
^{\infty}q^k\,R_n(q^k;q)\left(1+(-1)^\ep q^{(A/2-1)(n+2k)}\right)\\
\notag
&=q^{1 + \frac{A n}{4} + \frac{rn }{2} + \frac{rA n^2 }{2} -
        r^2 n^2 }
  \frac{
     {({\let \over / \def\frac#1#2{#1 / #2} q}; q) _{n} }^{A - 2 r}\,
     ({\let \over / \def\frac#1#2{#1 / #2} q}; q) _{r n} \,
     ({\let \over / \def\frac#1#2{#1 / #2} q^ {(r+1)n+2}}; q)
      _{r n}  }{{(\let \over /
       \def\frac#1#2{#1 / #2} q^ { rn+1}}; q) _{n+1} ^A}\\
\notag
&\times
\sum _{k=0} ^{\infty}\left(1+(-1)^\ep q^{(A/2-1)((2r+1)n+2k+2)}\right)
\frac {(q^ { (2r+1)n+2};q)_k\,(q^ { r
n+1};q)_k^{A+1}} {(q^ {(r+1)n+2};q)_k^{A+1}}
\left(q^ { \frac{1}{2}(A - 2 r)n+1}\right)^k.
\end{align}
{\`A} cause du facteur $(1+(-1)^\ep q^{(A/2-1)((2r+1)n+2k+2)})$, il n'y 
a pas d'{\'e}criture {\'e}l{\'e}gante de cette somme en notation
hyperg{\'e}om{\'e}trique basique car il faudrait pour cela utiliser des 
racines $(A-2)$-i{\`e}mes de l'unit{\'e}. 
La somme est cependant une s{\'e}rie
bien {\'e}quilibr{\'e}e (``well-poised"; voir \cite[Sec.~2.1]{gr}),
et m{\^e}me tr{\`e}s bien {\'e}quilibr{\'e}e (``very-well-poised") 
si $\ep$ est impair.
Remarquons que la s{\'e}rie alternative $\tilde S_n(q)$ donn{\'e}e par 
\eqref{eq:Salt} est une s{\'e}rie tr{\`e}s bien {\'e}quilibr{\'e}e. 

\medskip

\subsection{Estimation asymptotique de $S_{n}^{[\ep]}(q)$}
Dans ce paragraphe on d{\'e}montre le point~B) du
paragraphe~\ref{para:resume}.

\begin{lem}\label{lem:asympserie}
Soit $A\ge r/2$.
Pour $\ep\in\{0,1\}$ et $q$ tel que $\vert q\vert<1$, on a
$$
\lim_{n\to +\io}\frac{1}{n^2}\log\vert S_{n}^{[\ep]}(q)\vert 
= -\frac 12 r(A-2r) \log\vert 1/q\vert .
$$
\end{lem}
\begin{proof} Fixons $A$, $r$ et $q$. 
Notons $\rho_{k}(q)=q^kR_n(q^k;q)$ le sommande de la s{\'e}rie
$S_n(q)$~: clairement, $\rho_k(q)=0$ pour 
$k\in\{1, \ldots, rn\}$ et $\rho_k(q)\ne 0$ pour 
tout $k\ge rn+1$. Il est imm{\'e}diat que pour $k\ge rn+1$, 
$$
\frac{\rho_{k+1}(q)}{\rho_{k}(q)}=
q^{(A-2r)n/2+1}\cdot 
\frac{1-q^{1+k+n+rn}}{1-q^{k-rn}}\cdot 
\left( \frac{1-q^k}{1-q^{k+n+1}}\right)^{A+1}.
$$
Donc, puisque $\vert q\vert <1$ et $k\ge rn+1$,  
$$
\left\vert\frac{\rho_{k+1}(q)}{\rho_{k}(q)}\right\vert \le \vert q\vert ^{(A-2r)n/2}\cdot 
\left(\frac{1+\vert q\vert }{1-\vert q\vert }\right)^{A+2}< \frac{1}{3},
$$
pourvu que $n\gg 0$ (ce qui d{\'e}pend de $A$, $r$ et $q$ mais pas de
$k$). {\`A} cette condition, 
on a d'un c{\^o}t{\'e} 
\begin{equation*}
\vert S_n(q)\vert = 
\vert \rho_{rn+1}(q)\vert \cdot
\bigg\vert \sum_{k=rn+1}^{\io}\frac{\rho_{k}(q)}{\rho_{rn+1}(q)}\bigg
\vert 
< \vert \rho_{rn+1}(q)\vert \cdot \sum_{k=0}^{\io} 3^{-k} = \frac 3 2 \,\vert
\rho_{rn+1}(q)\vert,
\end{equation*}
et d'un autre c{\^o}t{\'e}, 
$$
\vert S_n(q)\vert\ge  |\rho_{rn+1}(q)| \left(1-\sum_{k=rn+2}^{\io} \left\vert  
\frac {\rho_k(q)} {\rho_{rn+1}(q)} \right\vert  \right)
\ge \frac 12
\vert \rho_{rn+1}(q)\vert ,
$$
Comme 
$$
\rho_{rn+1}(q)= (q;q)_{n}^{A-2r}q^{-(A-2r)n/4}
\frac{(q;q)_{rn}\,(q^{(r+1)n+2};q)_{rn}}{(q^{rn+1};q)_{n+1}^A} \,q^{(A-2r)(rn+1)n/2}
$$
et 
$$
\lim_{n\to+\io} \frac{1}{n^2}\log\left\vert (q;q)_{n}^{A-2r}q^{-(A-2r)n/4}
\frac{(q;q)_{rn}\,(q^{(r+1)n+2};q)_{rn}}{(q^{rn+1};q)_{n+1}^A}\right\vert =0,
$$
on obtient   
\begin{equation*}
\lim_{n\to +\io}\frac{1}{n^2}\log\vert S_{n}(q)\vert 
=
 \lim_{n\to +\io}\frac{1}{n^2}
\log\vert q^{(A-2r)(rn+1)n/2}\vert= -\frac12 r(A-2r)  \log\vert 1/q\vert. 
\end{equation*}

La s{\'e}rie qui nous int{\'e}resse {\`a} proprement parler est $S^{[\ep]}_n(q)$,
que l'on peut {\'e}crire, selon \eqref{eq:Sexplicit}, sous la forme
$$
S^{[\ep]}_n(q)=
\sum_{k=rn+1}^{\io}\bigl(1+(-1)^{\ep} q^{(A/2-1)(n+2k)}\bigr)\rho_k(q).
$$
Le facteur $1+(-1)^{\ep} q^{(A/2-1)(n+2k)}$  n'a aucune influence 
asymptotique par rapport {\`a} $\rho_k(q)$ quand $n\to +\io$ et est toujours 
non nul 
puisque $\vert q\vert<1$.
En proc{\'e}dant exactement comme ci-dessus, on montre alors que
 le terme d'indice $k=rn+1$ domine les suivants et que  
$$
\lim_{n\to +\io}\frac{1}{n^2}\log\vert S_{n}^{[\ep]}(q)
\vert = -\frac12 r(A-2r) \log\vert  1/q\vert.
$$
\end{proof}

\subsection{Estimation asymptotique de $P_{s,n}^{[\ep]}(q)$}
\label{sec:P}
Dans ce paragraphe on d{\'e}montre le point~C) du
paragraphe~\ref{para:resume}.

\begin{lem}\label{lem:asympoly}
Pour tous $\ep\in \{0,1\}$, $s\in\{0, \ldots, A\}$ et 
$q$ tel que $\vert q\vert<1$, on a
$$
\limsup_{n\to +\io}\frac{1}{n^2}\log\vert P_{s,n}^{[\ep]}(q)\vert \le
\frac18(A+4r^2)\log \vert 1/q\vert.
$$
\end{lem}
\begin{proof}
{\'E}tant donn{\'e}es les expressions des diverses fonctions $P_{s,n}^{[\ep]}(q)$ 
apparues au court de  la d{\'e}monstration du Lemme~\ref{lem:dichoto}, 
il suffit de montrer que l'estimation attendue est d{\'e}j{\`a} valide pour les coefficients
$d_{s,j,n}(q)=(-1)^sq^{js}c_{s,j,n}(q)$, uniform{\'e}ment en $j$ et $s$.  
Fixons $q$ tel que $\vert q\vert <1$ et $j\in\{0,\ldots,n\}$ et  
soit $\eta>0$ assez petit
(on peut choisir $\eta=(1-\vert q\vert)/2$)~: 
par la formule de Cauchy
appliqu{\'e}e {\`a}~\eqref{eq:c}, on obtient 
$$
c_{s,j,n}(q)
=\frac{1}{2\pi i}\int_{\mathcal{C}}R_n(T;q)(T-q^{-j})^{s-1}\dd T,
$$
o{\`u} $\mathcal{C}$ d{\'e}signe le cercle de centre $q^{-j}$ et de rayon
$\eta\vert q\vert^{-j}$. En rempla\c cant $T$ par $Tq^{-j}$, cela {\'e}quivaut {\`a}
$$
d_{s,j,n}(q)
=-\frac{1}{2\pi i}\int_{\mathcal{C'}}R_n(Tq^{-j};q)
(1-T)^{s-1}\dd T,
$$
o{\`u} $\mathcal{C'}$ 
d{\'e}signe le cercle de centre $1$ et de rayon
$\eta$. 
Il s'agit donc de majorer la fraction 
\begin{align}
\label{eq:Rn}
R_n(Tq^{-j}&;q)(1-T)^{s-1}\\
\notag
&= q^{-(A-2r)(2j+1)n/4}\,(q;q)_n^{A-2r}
\frac{(q^{-rn-j}T;q)_{rn}(q^{n-j+1}T;q)_{rn}\,(1-T)^{s-1}}
{(Tq^{-j};q)_{n+1}^{A}}\;
T^{(A-2r)n/2}\\
\notag
&={\left( -1 \right) }^{Aj+ n r} 
q^{A\binom {j+1}2- \frac{A n}{4} (2j+1) - \frac{ r^2n^2}{2}} 
 ({\let \over / \def\frac#1#2{#1 / #2} q}; q) _{n} ^{A - 2 r}\\
\notag
&\kern4cm
\times
 \frac{
     ({\let \over / \def\frac#1#2{#1 / #2} \frac{q^ {j+1}}{T}};
      q) _{rn} \,({\let \over / \def\frac#1#2{#1 / #2} q^ {n -
      j + 1} T}; q) _{rn}\,
     {\left( 1 - T \right) }^{ s-A-1}  }
{({\let \over / \def\frac#1#2{#1 / #2}
        \frac{q}{T}}; q) _{j} ^A\,
     ({\let \over / \def\frac#1#2{#1 / #2} q T}; q) _{n-j} ^A}
 T^{ \frac{A }{2}(n-2 j) } 
\end{align}
sur $\mathcal{C'}$. On v{\'e}rifie que pour tous nombres entiers positifs
$a,b$, avec $a>0$, et tout $T\in\mathcal{C}$, on a  
$$
0<(\vert q\vert(1+\eta);\vert q\vert)_\infty\le
\vert(q^aT;q)_b\vert\le (-(1+\eta);\vert q\vert)_\infty
$$
et
$$
0<(\vert q\vert/(1-\eta);\vert q\vert)_\infty\le
\vert(q^a/T;q)_b\vert\le (-1/(1-\eta);\vert q\vert)_\infty.
$$
Les bornes inf{\'e}rieures sont effectivement $>0$ parce que les
valeurs $\vert q\vert(1+\eta)$ et $\vert q\vert/(1-\eta)$ sont $<1$.
On a {\'e}galement  
$$
\vert (q;q)_n\vert\le (-\vert q\vert;\vert q\vert)_\infty\quad
\text{et}\quad 
\vert T^{(n-2j)A/2}\vert \le \bigl(\max\{1+\eta,1/(1-\eta)\}\bigr)^{An/2}.
$$
Finalement, pour l'exposant de $q$ dans le membre final de 
\eqref{eq:Rn}, on a
$$A\binom {j+1}2- \frac{A n}{4} (2j+1) - \frac{ r^2n^2}{2}\ge 
-\frac {A} {8}-\frac {An^2} {8}- \frac{ r^2n^2}{2},$$
pour tout $j\in\{0,1, \cdots, n\}$. 
On en d{\'e}duit que 
$$
\vert d_{s,j,n}(q)\vert \le 
C_0\cdot\vert q\vert^{-(A+4r^2)n^2/8},
$$
o{\`u} $C_0$ d{\'e}pend de $A$, $\eta$, $n$, $q$ et $r$ mais  ni de $j$ ni de $s$, et v{\'e}rifie 
$C_0^{1/n^2}\to 1$ quand $n\to +\io$. 
\end{proof}

\begin{remark}
Nous pensons qu'en fait 
$$\lim_{n\to+\infty}\frac{1}{n^2}\log\vert P_{s,n}^{[\ep]}(q)\vert =
\frac18(A+4r^2)\log \vert 1/q\vert,$$
ce qui est support{\'e} par des calculs num{\'e}riques effectu{\'e}s sur
ordinateur, ainsi que par des arguments
heuristiques. Si cela est vrai, on pourra alors {\'e}tendre les
Th{\'e}or{\`e}mes~\ref{theo:dimietp} et \ref{theo:modu} aux valeurs
alg{\'e}briques de $q$ en utilisant un crit{\`e}re d'ind{\'e}pendance
lin{\'e}aire sur les corps de nombres, d{\^u} {\`a} Jadot \cite{jadot}.
\end{remark}

\subsection{Estimation arithm{\'e}tique de $P_{s,n}^{[\ep]}(q)$} 
\label{ssec:arith}

Dans ce paragraphe on d{\'e}montre le point~D) du
paragraphe~\ref{para:resume}.
Laissons temporairement $q$ {\^e}tre  tel que $\vert q\vert \ne 1$.
Rappelons que l'on d{\'e}finit les coefficients $q$-binomiaux par 
$$
\begin{bmatrix} n\\  k\end{bmatrix}_{q}=\frac{(q;q)_n}{(q;q)_k\,(q;q)_{n-k}}
$$
et qu'ils appartiennent {\`a} $\mathbb{Z}[q]$.   
Pour tout $n\ge 0$, le polyn{\^o}me unitaire $d_n(q)$ de plus petit degr{\'e} 
en $q$ et tel que  
$$
d_n(q)\,\frac{1}{1-q},\; d_n(q)\,\frac{1}{1-q^2}, \ldots,\; 
d_n(q)\,\frac{1}{1-q^n} \in\mathbb{Z}[q],
$$ 
est donn{\'e} par 
$\dis d_n(q)=\prod_{\ell=1}^n \Phi_{\ell}(q)$, o{\`u} 
$\dis \Phi_{\ell}(q)=\prod_{k=1,\,(k,\ell)=1}^{\ell}(q-e^{2i\pi k/\ell})$ 
est le $\ell^{\;\textup{i{\`e}me}}$ polyn{\^o}me cyclotomique.   
Il satisfait l'estimation suivante pour tout $q$ tel que $\vert q\vert
>1$ (voir \cite[\S2]{buv} et \cite[Lemma~2]{vana})~: 
\begin{equation} \label{eq:dnasy} 
\lim_{n\to +\io} \frac{1}{n^2}\log\vert d_n(q)\vert=\frac{3}{\pi^2}\log\vert q\vert.
\end{equation}

Nous sommes maintenant en position d'expliciter 
les fractions rationnelles $D_n(q)$ qui
apparaissent au point~D) du paragraphe~\ref{para:resume} 
et d'{\'e}tablir leur comportement asymptotique lorsque $n$ tend vers
$+\infty$. Soit
\begin{equation} \label{eq:Dndef}
D_n(q)=(A-1)!\, q^{(A-2r)n/4-\lceil A (n+1)^2/8\rceil
-\frac {r^2n^2} {2}+\frac {rn} {2}- (A-1)n}\,
d_n(1/q)^{A}.
\end{equation}
On a alors le lemme suivant.

\begin{lem} \label{lem:arith}
Soit $A$ pair. 

\begin{itemize}
\item[\em i)] 
Pour tous $\ep\in\{0,1\}$, $s\in\{0,\ldots A\}$ et $q$ tel que $\vert
q\vert\ne 1$, on a
$$
D_n(q)P_{s,n}^{[\ep]}(q)\in\mathbb{Z}[1/q];
$$
\item[\em ii)] 
Lorsque $\vert q\vert <1$, on a alors
$$
\lim_{n\to +\io} \frac{1}{n^2}\log \vert D_n(q) \vert
= \biggl(\frac{3}{\pi^2}A+\frac A8 + \frac {r^2} 2\biggr)\log \vert 1/q\vert.
$$ 
\end{itemize}
\end{lem}

\begin{proof} Nous allons d'abord d{\'e}montrer que pour 
tous $s\in\{0,\ldots A\}$, $j\in\{0, \dots, n\}$ et $q$ tel que $\vert q\vert\ne 1$, on a
$$
q^{(A-2r)n/4-\lceil A (n+1)^2/8\rceil
-\frac {r^2n^2} {2}+\frac {rn} {2}- n}\,
d_n(1/q)^{A-s} c_{s,j,n}(q)\in\mathbb{Z}[1/q],
$$ 
le coefficient $c_{s,j,n}(q)$ {\'e}tant donn{\'e} par \eqref{eq:c}.
D{\'e}composons
  tout d'abord $R_n(T;q)$ de la fa{\c c}on suivante~:
\begin{multline*}
R_n(T;q)(T-q^{-j})^A=q^{-(A-2r)n/4}\left(q^{-n(n+1)/2}(q;q)_n
\frac{T^n\,(T-q^{-j})}{(T-1)\cdots(T-q^{-n})}
\right)^{(A-2r)/2}\\
\cdot
\left(q^{-n(n+1)/2}(q;q)_n\frac{(T-q^{-j})}{(T-1)\cdots(T-q^{-n})}
\right)^{(A-2r)/2}\\
\cdot\prod_{\ell=1}^r\left(q^{-n(n+1)/2}\frac{(q^{-\ell n}T;q)_n\,(T-q^{-j})}
{(T-1)\cdots(T-q^{-n})}\right)\\
\cdot\prod_{\ell=1}^r\left(q^{-n(n+1)/2}\frac{(q^{\ell n+1}T;q)_n\,(T-q^{-j})}
{(T-1)\cdots(T-q^{-n})}\right).
\end{multline*}
On v{\'e}rifie en d{\'e}composant en {\'e}l{\'e}ments simples que l'on a les
quatre identit{\'e}s suivantes~:
\begin{align} \label{eq:F1}
F(T)&=q^{-n(n+1)/2}(q;q)_n\frac{T^n\,(T-q^{-j})}{(T-1)\cdots(T-q^{-n})}\\
\notag
&=(-1)^n(1/q;1/q)_n+
\underset{i\ne j}{\sum_{{i=0}}^n} (-1)^{n+i} 
q^{-i(i+1)/2}\begin{bmatrix} n\\  i\end{bmatrix}_{1/q}
\frac{q^{-i}-q^{-j}}{T-q^{-i}}
\end{align}
\begin{align} 
G(T)=\frac{q^{-n(n+1)/2}\,(q;q)_n\,(T-q^{-j})}{(T-1)\cdots(T-q^{-n})}
\label{eq:G1}
&=
\underset{i\ne j}{\sum_{{i=0}}^n} (-1)^{n+i}
q^{in-i(i+1)/2}\begin{bmatrix} n\\  i\end{bmatrix}_{1/q}
\frac{q^{-i}-q^{-j}}{T-q^{-i}}
\end{align}
\begin{align}  \label{eq:H1}
H_{\ell}(T)&=q^{-n(n+1)/2}
\frac{(q^{-\ell n}T;q)_n\,(T-q^{-j})}
{(T-1)\cdots(T-q^{-n})}\\
\notag
&=(-1)^n q^{-\ell n^2-n}+\underset{i\ne j}{\sum_{{i=0}}^n}
(-1)^{i} 
q^{-(n-i)^2/2-(n+i)/2}
\begin{bmatrix} n\\  i\end{bmatrix}_{1/q}
\begin{bmatrix} \ell n+i\\ n\end{bmatrix}_{1/q}
\frac{q^{-i}-q^{-j}}{T-q^{-i}}
\end{align}
et
\begin{align}  \label{eq:I1}
I_{\ell}(T)&=q^{-n(n+1)/2}
\frac{(q^{\ell n+1}T;q)_n\,(T-q^{-j})}
{(T-1)\cdots(T-q^{-n})}\\
\notag
&=
(-1)^n q^{\ell n^2}\biggl(1+\underset{i\ne j}{\sum_{{i=0}}^n}
(-1)^{i} 
q^{-i(i+1)/2}
\begin{bmatrix} n\\  i\end{bmatrix}_{1/q}
\begin{bmatrix} \ell(n+1)-i\\ n\end{bmatrix}_{1/q}
\frac{q^{-i}-q^{-j}}{T-q^{-i}}\biggr)
\end{align}
Notons 
$\displaystyle \partial^{\mu }=
\frac{1}{\mu !}\frac{\dd^{\mu}}{\dd T^{\mu}}$.  Il r{\'e}sulte des formules 
\eqref{eq:F1}, \eqref{eq:G1}, \eqref{eq:H1}, \eqref{eq:I1},  
que pour tout entier 
$\mu\ge 0$, tout $j\in\{0,\ldots, n\}$ et tout $\ell\in\{1,\ldots, r\}$, les quantit{\'e}s 
\begin{equation*}
d_n(1/q)^{\mu}\partial^{\mu} F(q^{-j}),\;\; 
d_n(1/q)^{\mu}\partial^{\mu}G(q^{-j}),\;\;
d_n(1/q)^{\mu}\partial^{\mu}H_{\ell}(q^{-j}),\;\;
d_n(1/q)^{\mu}\partial^{\mu}I_{\ell}(q^{-j}),
\end{equation*}
sont des polyn{\^o}mes en $1/q$, multipli{\'e}s par une certaine puissance de
$q$ et dont 
les coefficients sont des nombres entiers. 
Or la formule de Leibniz nous
donne 
\begin{multline*}
c_{s,j,n}(q)=q^{-(A-2r)n/4}\sum_{\boldsymbol{\mu}}(\partial^{\mu_1}F)
\cdots(\partial ^{\mu_{(A-2r)/2}}F)
(\partial^{\mu_{(A-2r)/2+1}}G)\cdots (\partial^{\mu_{A-2r}}G)\\\cdot
(\partial^{\mu_{A-2r+1}}H_1)\cdots (\partial^{\mu_{A-r}}H_r)
(\partial^{\mu_{A-r+1}}I_1)\cdots (\partial^{\mu_{A}}I_r),
\end{multline*}
(o{\`u} la sommation est sur tous les multiplets
$\boldsymbol{\mu}=(\mu_1,\ldots,\mu_{A})$ tels que
$\mu_1+\cdots+\mu_{A}=A-s$, et o{\`u} l'on a omis l'{\'e}valuation en
$T=q^{-j}$ pour simplifier). On en d{\'e}duit que 
$$
q^{(A-2r)n/4}d_n(1/q)^{A-s}  c_{s,j,n}(q)
$$
est un polyn{\^o}me en $1/q$, {\`a} coefficients entiers et multipli{\'e} par une 
certaine puissance de $q$. 

Il nous faut maintenant 
expliciter cette puissance de $q$, c'est-{\`a}-dire d{\'e}terminer un exposant $e(n)$ le plus
grand possible tel que 
$$
q^{e(n)}q^{(A-2r)n/4}d_n(1/q)^{A-s}  c_{s,j,n}(q)\in \mathbb{Z}[1/q].
$$

Pour cela, revenons {\`a} la d{\'e}finition originelle \eqref{eq:c} de
$c_{s,j,n}(q)$. Posons pour simplifier $R(T)=R_n(T;q)(T-q^{-j})^A$.
En utilisant la formule de Fa{\`a} di Bruno
(voir \cite{jo}), on a d'un c{\^o}t{\'e} 
\begin{align} \label{eq:Rn1}
\partial^{\mu}R(T)&=\partial^{\mu}\exp(\log(R(T)))\\
\notag
&=
\sum _{} ^{}\frac {1} {k_1!\,k_2!\cdots k_\mu!}\left(\frac
{\partial^{\vert \mathbf k\vert}}
{\partial T^{\vert \mathbf k\vert}}\exp\right)\Bigl(\log(R(T))\Bigr)
\prod _{i=1} ^{\mu}\left({\partial^i
\log(R(T))} \right)^{k_i}\\
\notag
&=
\sum _{} ^{}\frac {1} {k_1!\,k_2!\cdots k_\mu!}R(T)
\prod _{i=1} ^{\mu}\left(\frac {1} {i!}\frac {\partial^{i-1}}
{\partial T^{i-1}} r(T) \right)^{k_i},
\end{align}
o{\`u} la somme porte sur tous les multiplets $\mathbf k=
(k_1,k_2,\dots,k_\mu)$ tels que
$k_1+2k_2+\dots+\mu k_\mu=\mu$, o{\`u} $\vert \mathbf k\vert$
d{\'e}signe la somme $k_1+k_2+\dots+k_\mu$, et o{\`u} $r(T)$ est la
d{\'e}riv{\'e}e logarithmique de $R(T)$,
\begin{align} \label{eq:r} 
r(T)=\frac{\frac {\partial}
{\partial T}R(T)} {R(T)}
&=\frac {(A-2r)n} {2}\frac {1} {T}+A\underset{i\ne j}{
\sum _{i=0} ^{n}}\frac {q^i} {1-q^iT}
-\sum _{i=-r n} ^{-1}\frac {q^i} 
{1-q^iT}-
\sum _{i=n+1} ^{(r+1)n}\frac {q^i} 
{1-q^iT}\\
\notag
&=
\frac {(A-2r)n} {2}\frac {1} {T}+A{
\sum _{i=0} ^{j-1}}\frac {q^i} {1-q^iT}
-A{
\sum _{i=j+1} ^{n}}\frac {1/T} {1-q^{-i}/T}\\
\notag
&\kern3cm
-\sum _{i=-r n} ^{-1}\frac {q^i} 
{1-q^iT}+
\sum _{i=n+1} ^{(r+1)n}\frac {1/T} 
{1-q^{-i}/T}.
\end{align} 
D'un autre c{\^o}t{\'e}, on a
\begin{equation} \label{eq:Rn2} 
R(T)\big\vert_{T=q^{-j}}=
{\left( -1 \right) }^{r n} 
 q^{- (A-2r)n/4-{A }\binom j2 +
        \frac{A j n}{2} -  \frac{ rn}{2}+ \frac{r^2n^2}{2}}
\frac{{
     ({  \frac{1}{q}}; \frac{1}{q})}
        _{n} ^{A - 2 r}\,{({\let \over / \def\frac#1#2{#1 / #2}
      q^{-j-1}}; \frac{1}{q})} _{rn} \,
    { ({\let \over / \def\frac#1#2{#1 / #2} q^{- n+j-1}};
      \frac{1}{q})} _{rn} }{{({
        \frac{1}{q}}; \frac{1}{q})} _{j} ^A\,
    { ({ \frac{1}{q}}; \frac{1}{q})}
        _{n-j} ^A}.
\end{equation}
Comme on sait d{\'e}j{\`a} que 
$$
q^{(A-2r)n/4}d_n(1/q)^{A-s}\partial^{A-s} R(T)\big\vert_{T=q^{-j}}
=q^{(A-2r)n/4}d_n(1/q)^{A-s}c_{s,j,n}(q) 
$$ 
est un polyn{\^o}me 
en $1/q$, multipli{\'e} par une puissance de $q$, et comme 
$$-\binom j2+\frac {jn} {2}\le \frac {1} {8}(n+1)^2,$$
la combinaison de \eqref{eq:Rn1}, \eqref{eq:r} et \eqref{eq:Rn2}
montre que
$$
q^{(A-2r)n/4-\lceil A (n+1)^2/8\rceil
-\frac {r^2n^2} {2}+\frac {rn} {2}- (A-s)n}\,
d_n(1/q)^{A-s}  c_{s,j,n}(q)
$$
est un polyn{\^o}me en $1/q$ {\`a} coefficients entiers. 
(Lorsque l'on utilise \eqref{eq:r}, 
des expressions comme $1/(1-q^{-m})^M$, o{\`u} $m$ est un nombre entier
positif, sont interpr{\'e}t{\'e}es comme s{\'e}ries formelles en $1/q$.)

On peut maintenant passer {\`a} la d{\'e}monstration du point~i)
proprement dit. De nouveau, on se concentre sur le cas $\ep=1$. 
Rappelons que pour $s\ge 2$ pair, on a $P_{s,n}(q)=0$ et que 
pour $s$ impair de $\{3, \ldots, A-1\}$, on a (voir \eqref{eq:S->P})
$$
\displaystyle P_{s,n}^{[1]}(q)=\sum_{k=s}^{A}\alpha_{k,s}P_{k,n}(1;q),
$$
avec 
$$
P_{k,n}(1;q)=\sum_{j=0}^n d_{k,j,n}(q)\,q^{-j} =
(-1)^k\sum_{j=0}^n 
q^{j(k-1)}c_{k,j,n}(q).
$$
Il r{\'e}sulte donc de l'{\'e}tude pr{\'e}c{\'e}dente sur $c_{s,j,n}(q)$ 
(et du d{\'e}nominateur commun $(A-1)!$ aux
$\al_{s,k}$) que pour $s\ge 3$~:
$$
(A-1)!\, q^{(A-2r)n/4-\lceil A (n+1)^2/8\rceil
-\frac {r^2n^2} {2}+\frac {rn} {2}- (A-s) n}\,
d_n(1/q)^{A-s} P_{s,n}^{[1]}(q)\in\mathbb{Z}[1/q].
$$
Par ailleurs,  
$$
\displaystyle P_{0,n}^{[1]}(q)=P_{0,n}(1;q)-
q^{-n}P_{0,n}(1;1/q)-\frac{\partial P_{1,n} }{\partial z}(1;q)
$$
avec~:
\begin{equation*}
P_{0,n}(1;q)
=-\sum_{s=1}^A\sum_{j=1}^n\sum_{k=1}^j
d_{s,j,n}(q)\frac {q^{k-j}} {(1-q^k)^s},
\end{equation*}
qui {\'e}quivaut {\`a}
\begin{equation*}
q^{-n}P_{0,n}(1;1/q)
=-\sum_{s=1}^A\sum_{j=0}^{n-1}\sum_{k=1}^{n-j}
(-1)^s c_{s,j,n}(q)\frac{q^{j(s-1)-k}}{(1-q^{-k})^s}
\end{equation*}
(en utilisant la relation $d_{s,j,n}(1/q)=d_{s,n-j,n}(q)$), et 
$$
\frac{\partial P_{1,n} }{\partial z}(1;q)=
-\sum _{j=0} ^{n}j\cdot c_{1,j,n}(q).
$$
On d{\'e}duit de ces trois expressions que 
$$
q^{(A-2r)n/4-\lceil A (n+1)^2/8\rceil
-\frac {r^2n^2} {2}+\frac {rn} {2}- (A-1)n}\,
d_n(1/q)^{A}P_{0,n}^{[1]}(q)\in\mathbb{Z}[1/q].
$$

Le point~ii) est {\'e}vident compte-tenu de~\eqref{eq:dnasy}.
\end{proof}

\section{$q$-Analogie avec la fonction z{\^e}ta de Riemann}
\label{sec:qanalogiezeta}

Nous terminons cet article en discutant plus en d{\'e}tail deux aspects 
du travail pr{\'e}sent{\'e} pr{\'e}c{\'e}demment. 
Dans la premi{\`e}re partie de ce paragraphe,  
nous justifions que nos $\zeta_q(s)$ sont des 
$q$-analogues normalis{\'e}s des nombres $\zeta(s)$, et que nos
s{\'e}ries $S^{[\ep]}_n(q)$ sont des $q$-analogues des s{\'e}ries   
utilis{\'e}es r{\'e}cemment dans les travaux sur la dimension
des espaces vectoriels engendr{\'e}s sur $\mathbb Q$ par les valeurs de
la fonction z{\^e}ta aux entiers impairs. Dans la deuxi{\`e}me partie, 
nous pr{\'e}sentons des $q$-analogues des s{\'e}ries, maintenant 
classiques, de Ball et Beukers--Gutnik--Nesterenko, utilis{\'e}es
dans des d{\'e}monstrations de l'irrationalit{\'e} de $\zeta(3)$~: 
ces $q$-analogues 
sont fortement li{\'e}s  {\`a} notre s{\'e}rie $S^{[1]}_n(q)$ pour $A=4$ 
et $r=1$.

\subsection{
{\og Convergence \fg} 
vers le cas {\og classique \fg}} \label{sec:class}
Dans ce paragraphe, on discute du rapport entre nos $q$-fonctions et
$q$-s{\'e}ries et les fonctions et s{\'e}ries {\og classiques \fg}.

Remarquons tout d'abord que les valeurs $\zeta_q(s)$ sont des 
$q$-analogues des valeurs de la fonction z{\^e}ta de Riemann
en nombres entiers $s$~: en utilisant les nombres de Stirling
de seconde esp{\`e}ce $S(N,j)$, d{\'e}finis par 
$$\ell^N=
\sum _{j=1} ^{N}S(N,j)\,\ell(\ell-1)\cdots(\ell-j+1),$$
on a
\begin{align*}
\zeta_q(s)&=
\sum _{k=1} ^{\infty}
\sum _{d=1} ^{\infty}d^{s-1}q^{k d}=\sum _{k=1} ^{\infty}
\sum _{d=0} ^{\infty}(d+1)^{s-1}q^{k (d+1)}\\
&=\sum _{k=1} ^{\infty}
\sum _{d=0} ^{\infty}
\sum _{j=1} ^{s-1}(-1)^{s-1-j}S(s-1,j)\,j!\binom {d+j}{j}q^{k
(d+1)}\\
&=\sum _{j=1} ^{s-1}(-1)^{s-1-j}S(s-1,j)\,j!\sum _{k=1} ^{\infty}
\frac {q^k} {(1-q^k)^{j+1}}.
\end{align*}
La s{\'e}rie $\mathcal{Z}_{j+1}(q)=\sum _{k=1} ^{\infty}
\frac {q^k} {(1-q^k)^{j+1}}$ {\'e}tant clairement un
$q$-analogue de $\zeta(j+1)$, la s{\'e}rie  $\zeta_q(s)$ est une
combinaison lin{\'e}aire de $q$-analogues des valeurs de la fonction
z{\^e}ta. On voit alors sur l'expression pr{\'e}c{\'e}dente que 
\begin{equation} \label{eq:Elim} 
\lim_{q\to1}(1-q)^s\,\zeta_q(s)=(s-1)!\,\zeta(s),
\end{equation}
et  $\zeta_q(s)$ peut donc 
aussi {\^e}tre consid{\'e}r{\'e}e comme un 
$q$-analogue normalis{\'e} de $\zeta(s)$ 
(voir \cite[Theorem~2]{jap} et \cite{zud5} pour d'autres 
d{\'e}monstrations de 
\eqref{eq:Elim}).

De m{\^e}me, en utilisant \eqref{eq:Z-Z->E}, \eqref{eq:Z+Z->E} et
\eqref{eq:Elim}, on obtient 
$$\lim_{q\to1}(1-q)^s\,(\mathcal Z_s(q)+(-1)^\ep\mathcal Z_s(1/q))
=2\zeta(s),$$
si $\ep\equiv s$ mod $2$, respectivement
$$\lim_{q\to1}(1-q)^{s-1}\,(\mathcal Z_s(q)+(-1)^\ep\mathcal Z_s(1/q))
=(s-2)\zeta(s-1),$$
si $\ep\not\equiv s$ mod $2$. Par cons{\'e}quent, notre combinaison
lin{\'e}aire $\mathcal Z_s(q)+(-1)^\ep\mathcal Z_s(1/q)$ (si
importante dans la d{\'e}monstration du Lemme~\ref{lem:dichoto}) est elle aussi 
un $q$-analogue normalis{\'e} de
$\zeta(s)$, respectivement $\zeta(s-1)$. En particulier, la combinaison
\begin{equation} \label{eq:zetaq3} 
\frac {1} {2}\bigl(\mathcal Z_3(q)-\mathcal Z_3(1/q)\bigr)=
\frac {1} {2}\zeta_q(3)=
\frac {1} {2}\sum _{k=1} ^{\infty}\frac {q^k(1+q^k)} {(1-q^k)^3}
\end{equation}
 r{\'e}appara{\^\i}tra au paragraphe suivant.

Comparons maintenant nos s{\'e}ries $S^{[\ep]}_n(q)$
avec les s{\'e}ries utilis{\'e}es pour la d{\'e}monstration du
Th{\'e}or{\`e}me~\ref{th:rivoal}. 
Si l'on multiplie $S^{[1]}_n(q)$ par $(1-q)^{A-1}$ et fait tendre $q$ 
vers 1, on obtient la s{\'e}rie
\begin{equation} \label{eq:serie}
\left({A} -{2}\right)\cdot n!^{A-2r} 
\sum _{k=1} ^{\infty}\left(k+\frac n2\right)
\frac {(k-rn)_{rn}\,(k+n+1)_{rn}} {(k)_{n+1}^A}
\end{equation}
o{\`u} on a encore utilis{\'e} les symboles de Pochhammer (voir la
d{\'e}monstration du Lemme~\ref{lem:crucial}). Aux termes $k+n/2$ et 
$({A} -{2})$ pr{\`e}s, il s'agit de la s{\'e}rie utilis{\'e}e dans \cite{ri, br} 
pour d{\'e}montrer le Th{\'e}or{\`e}me~\ref{th:rivoal}. 
Au terme $A-2$ pr{\`e}s, 
la s{\'e}rie \eqref{eq:serie} est aussi un cas sp{\'e}cial d'une
s{\'e}rie introduite dans \cite[p.~51]{ri2} (et ensuite g{\'e}n{\'e}ralis{\'e}e
dans \cite[Sec.~8]{zud0}).

De fa{\c c}on similaire, si 
l'on multiplie $S^{[0]}_n(q)$ par $(1-q)^{A}$ et fait tendre $q$ 
vers 1, on obtient la s{\'e}rie
\begin{equation*}
2\cdot n!^{A-2r}\sum _{k=1} ^{\infty}
\frac {(k-rn)_{rn}\,(k+n+1)_{rn}} {(k)_{n+1}^A},
\end{equation*}
qui est ici exactement celle introduite dans \cite{ri, br}.

Confront{\'e} {\`a} cette analogie frappante entre le {\og~cas $q$ g{\'e}n{\'e}ral ~\fg} 
et le {\og~cas $q=1$~\fg}, 
on aurait pu esp{\'e}rer analogie entre les minorations des dimensions des espaces vectoriels
engendr{\'e}s par les {\og~$q$-z{\^e}tas~\fg}, respectivement les
{\og~z{\^e}tas~\fg}, ce
qui n'est pas le cas (comparer les
Th{\'e}or{\`e}mes~\ref{th:rivoal} et \ref{theo:dimietp})~: 
les
m{\'e}thodes appliqu{\'e}es ne permettent pas de d{\'e}montrer des
minorations analogues.

\subsection{Autour de $\zeta_q(3)$} \label{sec:zetaq3}
On conna{\^\i}t actuellement de nombreuses d{\'e}monstrations de l'irrationalit{\'e}
de $\zeta(3)$ (voir l'article de synth{\`e}se \cite{fi}), dont plusieurs 
utilisent des s{\'e}ries hyperg{\'e}om{\'e}triques. 
Deux s{\'e}ries classiques sont la s{\'e}rie de Ball 
(qui est le cas $A=4$, $r=1$ 
de la s{\'e}rie \eqref{eq:serie}; voir l'introduction
de \cite{ri}),
\begin{equation} \label{eq:Ball} 
n!^2\sum_{k=1}^{\io}(2k+n)\frac{(k-n)_n\,(k+n+1)_n}{(k)_{n+1}^4},
\end{equation}
et la s{\'e}rie de Beukers--Gutnik--Nesterenko (voir 
\cite{beuk,gu,ne2}),
\begin{equation} \label{eq:bgn} 
-\sum_{k=1}^{\io}\frac{\dd }{\dd k}\left(\frac{(k-n)_n^2}{(k)_{n+1}^2}\right),
\end{equation}
qui sont en fait {\'e}gales pour tout nombre entier positif $n$.
Dans les preuves d'irrationalit{\'e} de $\zeta(3)$, on montre que ces
s{\'e}ries sont une combinaison lin{\'e}aire de $1$ et $\zeta(3)$, {\`a} 
coefficients rationnelles satisfaisant certaines propri{\'e}t{\'e}s asymptotiques et arithm{\'e}tiques. 
Dans ce paragraphe nous pr{\'e}sentons des $q$-analogues de ces deux 
s{\'e}ries, et nous d{\'e}montrons {\'e}galement que ces $q$-analogues sont 
une combinaison lin{\'e}aire de $1$ et de $\zeta_q(3)$. 

Nous commen\c cons avec la version 
non termin{\'e}e de la transformation de Watson due {\`a} Bailey
(voir \cite[(2.10.10); Appendix (III.36)]{gr})~:
\begin{multline} \label{eq:Bailey}
 {} _{8} \phi _{7} \! \left [                \begin{matrix} \let \over /
    \def\frac#1#2{#1 / #2} a, {\sqrt{a}} q, -{\sqrt{a}} q 
    , b, c, d, e, f\\ \let \over / \def\frac#1#2{#1 / #2} {\sqrt{a}},
    -{\sqrt{a}}, \frac{a q}{b}, \frac{a q}{c}, \frac{a q}{d},
    \frac{a q}{e}, \frac{a q}{f}\end{matrix} ;q, {\displaystyle
    \frac{a^ 2 q^ 2}{b c d e f}} \right ] \\=
     \frac {(\let \over / \def\frac#1#2{#1 / #2} a q, \frac{a q}{d e},
      \frac{a q}{d f}, \frac{a q}{e f} ;q) _\infty} {(\let \over /
      \def\frac#1#2{#1 / #2} \frac{a q}{d}, \frac{a q}{e},
      \frac{a q}{f}, \frac{a q}{d e f} ;q) _\infty} 
   {} _{4} \phi _{3} \! \left [                \begin{matrix} \let \over /
      \def\frac#1#2{#1 / #2} \frac{a q}{b c}, d, e, f\\ \let \over /
      \def\frac#1#2{#1 / #2} \frac{a q}{b}, \frac{a q}{c},
      \frac{d e f}{a}\end{matrix} ;q, {\displaystyle q} \right ] \\
 +
\frac {(\let \over / \def\frac#1#2{#1 / #2} a q,
      \frac{a q}{b c}, d, e, f, \frac{a^ 2 q^
      2}{b d e f}, \frac{a^ 2 q^ 2}{c d e f} ;q)
      _\infty} {(\let \over / \def\frac#1#2{#1 / #2} \frac{a q}{b},
      \frac{a q}{c}, \frac{a q}{d}, \frac{a q}{e}, \frac{a q}{f},
      \frac{a^ 2 q^ 2}{b c d e f}, \frac{d e f}{a q}
      ;q) _\infty}\\
\times
    {} _{4} \phi _{3} \! \left [                \begin{matrix} \let \over /
      \def\frac#1#2{#1 / #2} \frac{a q}{d e}, \frac{a q}{d f},
      \frac{a q}{e f}, \frac{a^ 2 q^ 2}{b c d e f}\\
      \let \over / \def\frac#1#2{#1 / #2} \frac{a^ 2 q^
      2}{b d e f}, \frac{a^ 2 q^ 2}{c d e f},
      \frac{a q^ 2}{d e f}\end{matrix} ;q, {\displaystyle q} \right
      ] .
\end{multline}
Dans cette identit{\'e}, on a utilis{\'e} la notation 
$$(a_1,a_2,\dots,a_k;q)_\infty:=(a_1;q)_\infty\,(a_2;q)_\infty\cdots
(a_k;q)_\infty.$$
On sp{\'e}cialise alors $a=\de^2q^{3n+2}$ et 
$b=c=d=e=f=\de q^{n+1}$ dans \eqref{eq:Bailey}, 
et on fait tendre $\de$ vers $1$. 
On obtient ainsi
\begin{multline*}
 \frac{({\let \over / \def\frac#1#2{#1 / #2} q^ {n+1}}; q)
         _{n+1} ^4 }{q^
       {{\left( n+1 \right) }^2}\,
      ({\let \over / \def\frac#1#2{#1 / #2} q}; q) _{n} \,
      ({\let \over / \def\frac#1#2{#1 / #2} q^ {2n+2}}; q)
       _{n+1} } \sum _{k=1} ^ {\infty} \frac{\left( 1 - q^ {2 k + n} \right)
       \,({\let \over / \def\frac#1#2{#1 / #2} q^ {k - n}}; q)
       _{n} \,({\let \over / \def\frac#1#2{#1 / #2} q^ {1 +
       k + n}}; q) _{n} }{({\let \over / \def\frac#1#2{#1 / #2}
       q^ k}; q) _{n+1} ^ 4}q^
       {k \left( n+1 \right) } \\
= \lim_{\de\to1}\Biggl(\frac{1
      }{q^{n+1}
        ({\let \over / \def\frac#1#2{#1 / #2} q}; q) _{n} }
 \frac {(\let \over / \def\frac#1#2{#1 / #2}
       {\de}^ 2 q^ {3 + 3 n}, q^ {n+1},
       q^ {n+1}, q^ {n+1}, \de q^ {
       n+1}, \de q^ {n+1}, \de q^ {
       n+1} ;q) _\infty} {(\let \over / \def\frac#1#2{#1 / #2}
       \de q^ {2n+2}, \de q^ {
       2 n+2}, \de q^ {2n+2},
       \de q^ {2n+2}, \de q^ {
       2 n+2}, \de q, \frac{q}{\de} ;q) _\infty}  \\
\times
\frac {1} {\de-1} \left( 
         \sum _{k=1} ^ {\infty} \frac{ ({\let \over /
          \def\frac#1#2{#1 / #2} q^ {k - n}}; q) _{n}
          \,({\let \over / \def\frac#1#2{#1 / #2}
          \de q^ {k - n}}; q) _{n} }{({\let \over /
          \def\frac#1#2{#1 / #2} \de q^ k}; q)
          _{n+1} ^ 2}(\de q^ k) -\sum _{k=1} ^ {\infty} \frac{\,({\let \over
         / \def\frac#1#2{#1 / #2} q^ {k - n}}; q) _{n}
         \,({\let \over / \def\frac#1#2{#1 / #2} \frac{q^ {k -
         n}}{\de}}; q) _{n} }{({\let \over /
         \def\frac#1#2{#1 / #2} q^ k}; q) _{n+1}
         ^ 2} q^ k \right)\Biggr),
\end{multline*}
soit, en utilisant le Th{\'e}or{\`e}me de l'H{\^o}pital, 
\begin{multline*}
 \frac{({\let \over / \def\frac#1#2{#1 / #2} q^ {n+1}}; q)
         _{n+1} ^4 }{q^
       {{\left( n+1 \right) }^2}\,
      ({\let \over / \def\frac#1#2{#1 / #2} q}; q) _{n} \,
      ({\let \over / \def\frac#1#2{#1 / #2} q^ {2n+2}}; q)
       _{n+1} } \sum _{k=1} ^ {\infty} \frac{\left( 1 - q^ {2 k + n} \right)
       \,({\let \over / \def\frac#1#2{#1 / #2} q^ {k - n}}; q)
       _{n} \,({\let \over / \def\frac#1#2{#1 / #2} q^ {1 +
       k + n}}; q) _{n} }{({\let \over / \def\frac#1#2{#1 / #2}
       q^ k}; q) _{n+1} ^ 4}q^
       {k \left( n+1 \right) } \\
= \frac{(q^{n+1};q)_n^4
      }{q^{n+1}
        ({\let \over / \def\frac#1#2{#1 / #2} q}; q) _{n} ^3\,
     (q^ {2n+2};q)_{n+1}}
         \sum _{k=1} ^ {\infty} \frac {\dd} {\dd\de}\left.\left(\frac{ 
          \,({\let \over / \def\frac#1#2{#1 / #2}
          \de q^ {k - n}}; q) _{n} ^2}{({\let \over /
          \def\frac#1#2{#1 / #2} \de q^ k}; q)
          _{n+1} ^ 2}(\de q^ k)  \right)\right\vert_{\de=1}.
\end{multline*}
Ceci {\'e}quivaut {\'e}videmment {\`a}
\begin{equation} \label{eq:id-Ball-bgn}
({\let \over / \def\frac#1#2{#1 / #2} q}; q) _{n}^2 
       \sum _{k=1} ^ {\infty} \frac{\left( 1 - q^ {2 k + n} \right)
       \,({\let \over / \def\frac#1#2{#1 / #2} q^ {k - n}}; q)
       _{n} \,({\let \over / \def\frac#1#2{#1 / #2} q^ {1 +
       k + n}}; q) _{n} }{({\let \over / \def\frac#1#2{#1 / #2}
       q^ k}; q) _{n+1} ^ 4}q^
       {k \left( n+1 \right) } 
= \frac{q^{n(n+1)}   }{\log q}
         \sum _{k=1} ^ {\infty} \frac {\dd} {\dd k}\left(\frac{ 
          \,({\let \over / \def\frac#1#2{#1 / #2}
           q^ {k - n}}; q) _{n} ^2}{({\let \over /
          \def\frac#1#2{#1 / #2}  q^ k}; q)
          _{n+1} ^ 2} q^ k  \right),
\end{equation}
dont le membre de gauche est un $q$-analogue parfait de la s{\'e}rie
\eqref{eq:Ball} de Ball et le membre de droite est un $q$-analogue parfait de la s{\'e}rie 
\eqref{eq:bgn} de Beukers--Gutnik--Nesterenko.
Plus pr{\'e}cis{\'e}ment, 
si l'on multiplie les deux membres de \eqref{eq:id-Ball-bgn} par
$(1-q)^3$ et que l'on fait tendre $q$ vers $1$, on retrouve l'{\'e}galit{\'e}
\eqref{eq:Ball}=\eqref{eq:bgn}. Il est {\'e}galement notable 
que, {\`a} un facteur $q^{-n/2}$ pr{\`e}s,   
le membre de gauche est exactement la s{\'e}rie 
$
S^{[1]}_n(q)=
S_n(q)-q^{-n}S_n(1/q)$ {\'e}tudi{\'e}e au 
paragraphe \ref{sec:demodimietp}, avec $A=4$ et $r=1$. 

Poursuivons l'analogie un peu plus loin en utilisant l'identit{\'e} 
int{\'e}grale de Agarwal (voir 
\cite[(4.6.5)]{gr}), 
\begin{multline*} 
 {} _{8} \phi _{7} \! \left [                \begin{matrix} \let \over /
    \def\frac#1#2{#1 / #2} q^a, q^{1+a/2}, -q^{1+a/2} 
    , q^b,q^ c, q^d, q^e, q^f\\ \let \over / \def\frac#1#2{#1 / #2}q^{a/2},
    -q^{a/2}, q^{1+a-b}, q^{1+a-c}, q^{1+a-d},
    q^{1+a-e},q^{1+a-f}\end{matrix} ;q, {\displaystyle
    q^{2+2a-b-c-d-e-f}} \right ] \\=
\sin\pi(b+c+d-a)
\frac
{(q^{1+a},q^b,q^c,q^d,q^{1+a-b-c},q^{1+a-b-d},q^{1+a-c-d},
 q^{1+a-e-f};q)_\infty} {(q,q^{b+c+d-a},q^{1+a-b-c-d},q^{1+a-b},
 q^{1+a-c},q^{1+a-d},q^{1+a-e},q^{1+a-f};q)_\infty}\\
\times\frac {1} {2\pi i}\int _{-i\infty} ^{i\infty}
\frac {(q^{1+s},q^{1+b-e+s},q^{1+a-f+s},q^{b+c+d-a+s};q)_\infty} 
{(q^{b+s},q^{c+s},q^{d+s},q^{1+a-e-f+s};q)_\infty}
\frac {\pi q^s\,ds} {\sin\pi s\,\sin\pi(b+c+d-a+s)},
\end{multline*}
o{\`u} la courbe d'int{\'e}gration est d{\'e}form{\'e}e autour de l'origine de sorte que
celle-ci soit {\`a} droite de la courbe. Cette identit{\'e} est  
 valable {\`a} condition que
$$\text{Re}\bigl(s\log q-\log(\sin\pi s\,\sin\pi(b+c+d-a+s))\bigr)<0,$$ 
et en  la sp{\'e}cia\-li\-sant en  
$a=3n+2$ et $b=c=d=e=f=n+1$, on peut donner une repr{\'e}sentation 
int{\'e}grale des
s{\'e}ries \eqref{eq:id-Ball-bgn}. 
Il s'agit en fait d'un cas limite 
puisqu'appara{\^\i}t le quotient 
$(\sin\pi(b+c+d-a))/(q^{1+a-b-c-d};q)_{\io}$~:
\begin{align*}
\underset{b,c,d\to n+1}{\lim_{a\to 3n+2}}\frac {\sin\pi(b+c+d-a)} 
{(q^{1+a-b-c-d};q)_{\io}}&=\lim_{a\to 3n+2}\frac {\sin\pi(3n+3-a)} 
{(q^{a-3n-2};q)_{\io}}\\
&=\lim_{x\to 0}\frac{\sin(\pi(1-x))}
{(1-q^x)(q^{1+x};q)_{\io}}=-\frac{\pi}{(q;q)_{\io}\log q}.
\end{align*}
On obtient ainsi~:
\begin{multline*}
{}_8 \phi_7 \left[
      \begin{array}{cccccccc}
      q^{3n+2} ,& q^{2+3n/2},& -q^{2+3n/2},& q^{n+1}, & q^{n+1}, & q^{n+1},& q^{n+1}, & q^{n+1}\\
                & q^{1+3n/2},& -q^{1+3n/2},& q^{2n+2},& q^{2n+2},& q^{2n+2},& q^{2n+2},& q^{2n+2} 
      \end{array}\; ;\,  q , q^{n+1} \right ]\\
=
\frac {(q^{n+1};q)_{n+1}^4} 
{(q;q)_n^3\,(q^{2n+2};q)_{n+1}}
\cdot \frac{1}{\log q}\cdot\frac{1}{2\pi i} \int\limits_{-i\io}^{i\io}
\frac
{(q^{1+s},q^{1+s},q^{2+2n+s},q^{2+2n+s}; q)_{\io}}
{(q^{1+n+s}, q^{1+n+s}, q^{1+n+s}, q^{1+n+s}; q)_{\io}}
\left(\frac{\pi}{\sin\pi s}\right)^2 q^s \dd s.
\label{eq:agarw}
\end{multline*}
 Par cons{\'e}quent
les s{\'e}ries \eqref{eq:id-Ball-bgn} s'expriment comme
\begin{equation} \label{eq:int}
\frac{q^{(n+1)^2}}{\log q}\cdot
\frac{1}{2\pi i} \int\limits_{-i\io}^{i\io}
\frac
{(q^{1+s},q^{1+s},q^{2+2n+s},q^{2+2n+s}; q)_{\io}}
{(q^{1+n+s}, q^{1+n+s}, q^{1+n+s}, q^{1+n+s}; q)_{\io}}
\left(\frac{\pi}{\sin\pi s}\right)^2 q^s \dd s.
\end{equation}
Remarquons que l'{\'e}galit{\'e} de \eqref{eq:int} avec le membre de
droite de \eqref{eq:id-Ball-bgn} peut {\^e}tre vue directement
en poussant la courbe d'int{\'e}gration {\`a} droite et 
en appliquant le th{\'e}or{\`e}me des r{\'e}sidus.  
On obtient ainsi un $q$-analogue parfait 
d'une autre identit{\'e} classique (voir \cite{ne2})~:
\begin{equation*} 
-\sum_{k=1}^{\io}\frac{\dd }{\dd k}\left(\frac{(k-n)_n^2}{(k)_{n+1}^2}\right)
= \frac{1}{2i \pi}\int\limits_{-i\io}^{i\io}
\frac
{\Gamma(1+s)^2\Gamma(2+2n+s)^2}
{\Gamma(1+n+s)^4}
\left(\frac{\pi}{\sin\pi s}\right)^2 \dd s.
\end{equation*}

Montrons maintenant que ces deux s{\'e}ries basiques peuvent s'{\'e}crire 
comme combinaison lin{\'e}aire de 1 et de $\zeta_q(3)$. Notons 
 $S_n(q)$ la s{\'e}rie {\`a} droite de \eqref{eq:id-Ball-bgn},
et 
$$
R_n(T ; q)= \frac{(q^{-n}T ; q)_n^2}{(T ; q)_{n+1}^2}
$$
(sans risque de confusion avec les notations 
utilis{\'e}es au paragraphe \ref{sec:demodimietp}) 
de sorte que 
$$
S_n(q)=
\sum_{k=1}^{\io}\frac{\dd }{\dd k}\left(q^kR_n(q^k ; q)\right).
$$
D{\'e}composons $R_n(T ; q)$ en {\'e}l{\'e}ments simples~: 
$$
R_n(T ; q)=\sum_{j=0}^n\left(\frac{a_{j,n}(q)}{(1-q^jT)^2}+\frac{b_{j,n}(q)}{1-q^jT}\right),
$$
d'o{\`u} 
$$
S_n(q)=\sum_{j=0}^n a_{j,n}(q) q^{-j}\sum_{k=1}^{\io} \frac{\dd }{\dd
    k}\left(\frac{q^{k+j}}{(1-q^{k+j})^2}\right)+
\sum_{j=0}^n b_{j,n}(q) q^{-j}\sum_{k=1}^{\io} \frac{\dd }{\dd
    k}\left(\frac{q^{k+j}}{1-q^{k+j}}\right).
$$
On a donc
\begin{multline*} 
S_n(q)=\left(\sum_{j=0}^n a_{j,n}(q) q^{-j}\right)\sum_{k=1}^{\io} \frac{\dd }{\dd
    k}\left(\frac{q^{k}}{(1-q^{k})^2}\right)+
\left(\sum_{j=0}^n b_{j,n}(q) q^{-j}\right)\sum_{k=1}^{\io} 
\frac{\dd }{\dd
    k}\left(\frac{q^{k}}{1-q^{k}}\right)\\ -\sum_{j=1}^n 
a_{j,n}(q) q^{-j}\sum_{k=1}^{j} \frac{\dd }{\dd
    k}\left(\frac{q^{k}}{(1-q^{k})^2}\right)-
\sum_{j=1}^n b_{j,n}(q) q^{-j}\sum_{k=1}^{j} \frac{\dd }{\dd
    k}\left(\frac{q^{k}}{1-q^{k}}\right).
\end{multline*}
Remarquons que 
$$
\sum_{j=0}^n b_{j,n}(q) q^{-j}=-\sum_{j=0}^n \textup{Res}(R_n( T;
q))_{T=q^{-j}}=\textup{Res}(R_n(T ; q))_{T=\io}=0
$$
et que l'on a
$$
\sum_{k=1}^{\io} \frac{\dd }{\dd
    k}\left(\frac{q^{k}}{(1-q^{k})^2}\right)= (\log q)\sum_{k=1}^{\io} 
\frac{q^{k}(1+q^k)}{(1-q^{k})^3}=(\log q)\zeta_q(3).
$$  
Finalement, en posant 
$$\displaystyle A_n(q)=\sum_{j=0}^n a_{j,n}(q) q^{-j}$$ 
et 
$$ 
B_n(q)=\sum_{j=1}^n\sum_{k=1}^{j} a_{j,n}(q) \frac{q^{k-j}(1+q^k)}{(1-q^{k})^3}+
\sum_{j=1}^n\sum_{k=1}^{j} b_{j,n}(q) \frac{q^{k-j}}{(1-q^{k})^2}, 
$$
on a donc 
$$
S_n(q)=(\log q)\bigl(A_n(q)\,\zeta_q(3)- B_n(q)\bigr), 
$$
c'est-{\`a}-dire 
\begin{equation}\label{eq:s=zetaq3}
\frac{1}{\log q}
         \sum _{k=1} ^ {\infty} \frac {\dd} {\dd k}\left(\frac{ 
          \,({\let \over / \def\frac#1#2{#1 / #2}
           q^ {k - n}}; q) _{n} ^2}{({\let \over /
          \def\frac#1#2{#1 / #2}  q^ k}; q)
          _{n+1} ^ 2} q^ k  \right)= A_n(q)\,\zeta_q(3)- B_n(q).
\end{equation}
Avec la technique expos{\'e}e au paragraphe~\ref{ssec:arith}, 
on montre qu'il existe $\bar D_n(q)\in\mathbb{Q}(q)$ tel que  
$\bar D_n(q)A_n(q)$ et $\bar D_n(q)B_n(q)$ soient dans $\mathbb{Z}[1/q]$ et 
$$
\lim_{n\to+\io}\frac{1}{n^2}\log\vert \bar D_n(q)\vert = \frac{9}{\pi^2}.
$$
Malheureusement, 
$$
\lim_{n\to+\io}\frac{1}{n^2}\log\vert S_n(q)\vert=0
$$
et on ne peut donc montrer l'irrationalit{\'e} de 
$\zeta_q(3)$ pour aucun $q\not =\pm1$ tel que $1/q\in\mathbb{Z}$. 

\begin{remark}
Pour finir, notons que, compte-tenu de~\eqref{eq:id-Ball-bgn}, 
on aurait {\'e}videmment pu utiliser le Lemme~\ref{lem:dichoto} 
avec $\ep=1$, $A=4$ et $r=1$ pour obtenir 
\eqref{eq:s=zetaq3} plus rapidement, mais au prix d'un $\bar D_n(q)$ 
moins bon~: cela 
sugg{\`e}re une probable {\og conjecture des d{\'e}nominateurs \fg}, 
comparable {\`a} celles 
{\'e}nonc{\'e}es dans le cas classique dans \cite{ri2} et \cite{zud0}.
Plus pr{\'e}cis{\'e}ment et comme not{\'e} dans la 
Remarque~(2) au paragraphe~\ref{para:resume}, on peut d{\'e}montrer les 
propri{\'e}t{\'e}s A)--D) du paragraphe~\ref{para:resume} pour
la s{\'e}rie alternative $\tilde S_n(q)$, en lieu et place de $S_n^{[\ep]}(q)$. 
Il semble alors que l'on peut dans ce cas remplacer le {\og 
d{\'e}nominateur\fg} $D_n(q)$ par un $\tilde D_n(q)$,
d{\'e}fini comme $D_n(q)$ par \eqref{eq:Dndef}, sauf que la puissance de
$d_n(q)$ est $A-1$ {\`a} la place de $A$. 
Comme pour $A=4$ et $r=1$, la 
s{\'e}rie $\tilde S_n(q)$ est aussi identique avec le membre gauche de 
\eqref{eq:id-Ball-bgn}, cette {\og conjecture des d{\'e}nominateurs\fg}
est v{\'e}rifi{\'e}e dans ce cas. Si elle l'est aussi 
pour $A=10$ et $r=2$, on pourra
d{\'e}montrer qu'au moins un des nombres $\zeta_q(3),\zeta_q(5),
\zeta_q(7),\zeta_q(9)$ est irrationnel, et donc am{\'e}liorer le
r{\'e}sultat du Th{\'e}or{\`e}me~\ref{theo:357911}.
Bien que pour $A=4$ et $r=1$ la s{\'e}rie $S_n^{[1]}(q)$ co{\"\i}ncide
avec $\tilde S_n(q)$ ({\`a} un facteur n{\'e}gligeable pr{\`e}s),
une telle am{\'e}lioration des d{\'e}nominateurs 
n'a probablement pas lieu, en g{\'e}n{\'e}ral,
pour la s{\'e}rie $S_n^{[\ep]}(q)$.
\end{remark}

\end{document}